\documentclass[a4paper,11pt]{article}

\usepackage[latin1]{inputenc}
\usepackage[T1]{fontenc}

\usepackage[english]{babel}
\usepackage{amsmath}
\usepackage{amsthm}
\usepackage{amssymb}
\usepackage{geometry}
\usepackage{mathtools}
\usepackage[draft,footnote,silent,marginclue]{fixme}
\usepackage{graphicx}
\usepackage{enumitem}
\usepackage{tikz}
\usetikzlibrary{calc} 

\theoremstyle{plain}
\newtheorem{thm}{Theorem}
\newtheorem{cor}{Corollary}
\newtheorem{prop}{Proposition}
\newtheorem{lem}{Lemma}
\newtheorem{fact}{Fact}

\theoremstyle{definition}
\newtheorem{defi}{Definition}

\theoremstyle{remark}
\newtheorem{rem}{Remark}

\newcommand{\pr}{\mathbb{P}}
\newcommand{\E}{\mathbb{E}}
\renewcommand{\d}[1]{\,\mathrm{d}#1}
\newcommand{\dd}[2]{\frac{\partial #1}{\partial #2}}
\newcommand{\ddtwo}[2]{\frac{\partial^2 #1}{\partial #2^2}}
\newcommand{\ddd}[3]{\frac{\partial^2 #1}{\partial #2\partial #3}}
\newcommand{\axisplus}{1}
\newcommand{\rhovalue}{}
\newcommand{\epsvalue}{}
\newcommand{\cvalue}{}
\renewcommand{\epsilon}{\varepsilon}

\newcommand{\limitshape}[3]{(1/(#1)*ln(exp(#1*#3)-exp(-#1*#3)-exp(-#1*(1-#2-#3))+
  exp(-#1*(#3-#2))) - 1/(#1)*ln(exp(#1*#2)-exp(-#1*(1-#2))))}

\DeclarePairedDelimiter{\paren}{(}{)}
\DeclarePairedDelimiter{\brac}{[}{]}

\DeclarePairedDelimiter{\abs}{\lvert}{\rvert}

\DeclarePairedDelimiter{\floor}{\lfloor}{\rfloor}

\author{Dan Beltoft\thanks{Aarhus University, Department of Math. Sciences, Ny Munkegade 118, 8000 Aarhus C, Denmark
},
C\'edric Boutillier\thanks{Universit\'e Pierre et Marie Curie, Laboratoire LPMA, 4 Place Jussieu 75005 Paris, France}$^{\phantom{\dag},}$\thanks{\'Ecole Normale Sup\'erieure, DMA, 45 rue d'Ulm 75005 Paris, France},
Nathana\"el Enriquez$^{\dag,}$\thanks{Universit\'e Paris-Ouest, Laboratoire MODAL'X, 200 Av. de la R\'epublique 92001 Nanterre, France}
}

\title{Random Young Diagrams in a Rectangular Box}

\begin{document}

\maketitle
\begin{abstract}
We exhibit the limit shape of random Young diagrams having a distribution proportional to the exponential of their area (grand-canonical ensemble), and confined in a rectangular box. The Ornstein-Uhlenbeck bridge arises from the fluctuations around the limit shape. The fluctuations for the unconfined case lead to a two-sided stationary Ornstein-Uhlenbeck process.
\end{abstract}



\section*{Introduction}

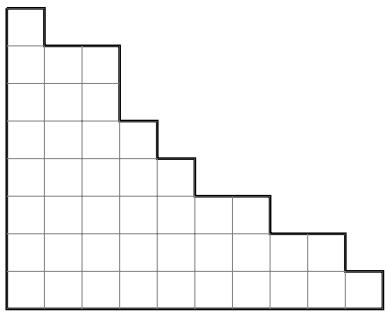
\begin{figure}[htb]
  \centering
  \begin{tikzpicture}[line width=1pt,scale=0.5]
    \draw[clip] (0,8) -| (1,7) -| (3,5) -| (4,4) -| (5,3) -| (7,2) -| (9,1) -| (10,0) -|
    (0,8); 
    \draw[help lines] (0,0) grid (10,8);
  \end{tikzpicture}
  \caption{The Young diagram for the partition $\pi=(8,7,7,5,4,3,3,2,2,1)$ of $n=42$.} 
  \label{fig:youngdiag}
\end{figure}

A \emph{partition} of an integer $n$ is a finite non-increasing sequence of integers
\begin{equation*}
  \pi=(\pi_1,\pi_2,\dots, \pi_k),
\end{equation*}
with $\pi_1+\cdots+\pi_k=n$. A convenient graphical representation in $(\mathbb{Z}^+)^2$
of the partition $\pi$ is its \emph{Young diagram}, also denoted by $\pi$. The Young
diagram consists of a stack of $\pi_i$ unit squares on the $i$-th column. See Figure
\ref{fig:youngdiag}. The \emph{area} $|\pi|$ of the Young diagram $\pi$ is the total
number of squares, which is equal to $n$.

If the number of summands $k$ is less than some integer $a$ and all the summands $\pi_j$
are less than $b$, then the Young diagram stays in the box $[0,a]\times [0,b]$.

In this paper we study the limit shape and fluctuations of large random Young diagrams
assigned to stay inside a (large) rectangular box, when the probability measure is
proportional to $q^{|\pi|}$, and $q$ is a suitable parameter. It is easy to see that if
one doesn't want the system to degenerate in the limit, one has to make $q$ depend on the
size of the box. Namely $1-q$ has to be taken of the order of the inverse of the side of
the box.

The study of the combinatorics of partitions of integers goes back to Hardy and Ramanujan
and then Erd\H{o}s in the 40's. The statistical physics point of view was introduced by
Vershik \cite{Vershik} in his study of the typical shape of the partition of a large
integer and obtained what we will call Vershik's curve. Recently, Funaki and Sasada obtained Vershik's result, in \cite{FunakiSasada}, as a by-product of an hydrodynamic result  for the corresponding particle system model.

Our study lies at the intersection of classical topics of probability theory,
combinatorics, and statistical physics.

Along the paper, we make an extensive use of the classical Gauss polynomials, well-known
in combinatorics as the generating functions for the number of Young diagrams with given
area (see e.g. \cite{AndrewsEriksson}). In order to get to the asymptotic regime, we are
led to state a $q$-analogue for Stirling's formula.

For all values of the parameters of the problem (limiting aspect ratio of the box, and
$c=\lim n(1-q)$), the limit shape obtained turns out to be a restriction of Vershik's
curve as it was recently noticed by Petrov in \cite{Petrov}.

The main part of the paper deals with the fluctuations around this limit shape and requires a fine understanding of the boundary of the random diagram. From a probabilistic point of view, the boundary is described by a Markov  which can be considered as the
$q$-analogue of the bridge of the simple random walk. In the classical case ($c=0$), the
properly rescaled interface converges to the Brownian bridge. In the general case, the
limiting process is an Ornstein-Uhlenbeck bridge.
Let us mention that, in the framework of unlimited partitions, and under some specifications on the
summands, Vershik and Yakubovich \cite{VershikYakubovich,Yakubovich} already observed Gaussian fluctuations. For strict partitions, see \cite{FreimanVershikYakub}.
The fluctuations for the unconfined case lead to a two-sided stationary Ornstein-Uhlenbeck process.
However, in contrast with the limit shape problem, there is no direct argument to deduce fluctuations for partitions in a box from the unlimited case.

The paper is organized as follows. After a presentation of the combinatorics of the
problem (Section 2), we perform its asymptotic analysis through a $q$-version of
Stirling's formula (Section 3). Section 4 is devoted to the derivation of the limit shape
phenomenon. In Section 5, we study the fluctuations around the limit shape: we first
compute the limiting 2-correlation function and deduce the convergence of
finite-dimensional marginals. This section ends with the delicate proof of the tightness
of the fluctuations.

\section{Presentation of the model}
\label{sec:present}

We study the asymptotic distribution of random Young diagrams fitting in a large
rectangular box with dimensions $a\times b$. Given a real number $q>0$, we
assign to each diagram $\pi$ the probability
\begin{equation*}
  \mathbb{P}^q_{a,b} (\pi) = \frac{q^{\abs{\pi}}}{Z_{a,b}(q)},
\end{equation*}
where $\abs{\pi}$ is the number of boxes of $\pi$ and $Z_{a,b}(q)$ is the
\emph{partition function}, the sum of all $q^{\abs{\pi}}$. For the sake of clarity, the probability
measure $\pr_{a,b}^q$ will be simply denoted, by $\pr$.

We fix a parameter $\rho\in(0,1)$, and choose the dimensions of the box for each
$n$ to be $a_n\times b_n$, where $(a_n)$ and $(b_n)$ are sequences of positive
integers satisfying:
\begin{equation*}
  a_n+b_n = 2n, \quad \lim_{n\rightarrow \infty}\frac{a_n}{b_n} = \frac{\rho}{1-\rho}.
\end{equation*}
We are interested in the limiting behavior when $n$ goes to $\infty$. To obtain
a non-degenerate limit, $q$ must go to $1$ as $n$ goes to infinity. We fix a
real parameter $c$ and pose $q=e^{-\frac{c}{n}}$. 

Note that $c=0$ corresponds to a uniform probability
measure. We will assume that $c\neq 0$. The results for $c=0$ can be obtained by taking limits. The physical meaning of the parameter $c$ is that of a pressure, since it is the variable conjugated to the (two-dimensional) volume.

\section{Combinatorics of the partitions}

Let us start with a fixed box with dimensions $a\times b$.
The partition function $Z_{a,b}(q)$ is expressed in terms of Gaussian
polynomials or $q$-binomial coefficients, where integers $j$ are replaced by their $q$-analogues $(j)_q=\frac{1-q^j}{1-q}$.
\begin{defi}
  Let $q>0$,
  \begin{gather*}
    n!_q = \prod_{j=1}^{n} (j)_q = \prod_{j=1}^{n}\frac{1-q^j}{1-q}\\
    \binom{n}{m}_q = \frac{n!_q}{(n-m)!_q m!_q} = \prod_{j=0}^{m-1}
    \frac{1-q^{n-j}}{1-q^{m-j}}.
  \end{gather*}
\end{defi}

\begin{lem}
  For all $(a,b)\in \paren*{\mathbb{N}^*}^2$, and all $q>0$, the partition
  function $Z_{a,b}(q)$ is equal to
  \begin{equation*}
    Z_{a,b}(q) =\binom{a+b}{a}_q = \binom{a+b}{b}_q.
  \end{equation*}
\end{lem}

\begin{proof}
  $\binom{a+b}{a}_q$ and $Z_{a,b}(q)$ both follow the recursion relation
  \begin{equation*}
    Z_{a,b}(q) = Z_{a,b-1}(q) + q^b Z_{a-1,b}(q),
  \end{equation*}
  where the first term corresponds to those diagrams with all parts strictly
  smaller than $b$, and the second term to the diagrams with at least one part
  of size $b$.
\end{proof}
   
We use the following coordinates. The bounding rectangle is the rectangle in the
plane with corners at the points $(0,0)$ and $(a+b,b-a)$ 
and sides with slopes $\pm 1$. The boundary of a diagram is encoded as a lattice path
$\paren[\big]{X_k}_{0\leq k\leq n}$ from the origin to the point $(a+b,b-a)$
and such that $X_{k+1}=X_k\pm 1$ for all $k$.
  
\begin{figure}[htb]
  \label{fig:latticepath}
  \centering
  \begin{tikzpicture}[line width=1pt,scale=.6]
    \renewcommand{\axisplus}{0.8}
     \begin{scope}[rotate=45]
       \draw[clip] (0,0) -- ++(1,0) -- ++(0,-2) -- ++(1,0) --
        ++(0,-2) coordinate (Xk) -- ++(4,0) -- ++(0,-1) -- ++(4,0) --
        ++(0,-1) -- ++(-10,0) -- ++(0,6); 
        \draw [style=help lines] (0,-6) grid (10,0);
      \end{scope}
      \begin{scope}[rotate=45]
        \draw (0,0) rectangle (10,-6) coordinate (rightcorner); 
        \path (10,0) coordinate (uppercorner); 
        \path (0,-6) coordinate (lowercorner);
      \end{scope}
      \path  (0,0 |- uppercorner) ++(0,\axisplus) coordinate (upperaxisend);
      \draw[->] (0,0 |- lowercorner) ++(0,-\axisplus) -- (upperaxisend); 
      \draw[<-] (0,0 -| rightcorner) ++(\axisplus,0) -- (-\axisplus,0);
      \draw (0,0 -| rightcorner) ++(0,-.2) node[below] {$2n$} -- ++(0,.4);
      \draw (0,0 |- lowercorner) ++(-.2,0) node[left] {$-a_n$} -- ++(.4,0);
      \draw (0,0 |- uppercorner) ++(-.2,0) node[left] {$b_n$} -- ++ (.4,0); 
      \fill (Xk) circle (.2);
      \draw (0,0 -| Xk) ++(0,-.2) node[below] {$2k$} -- ++(0,.4);
      \draw (0,0 |- Xk) ++ (-.2,0) node[left] {$X_{2k}$} -- ++(.4,0);
  \end{tikzpicture}
  \caption{A partition as a lattice path.}
\end{figure}
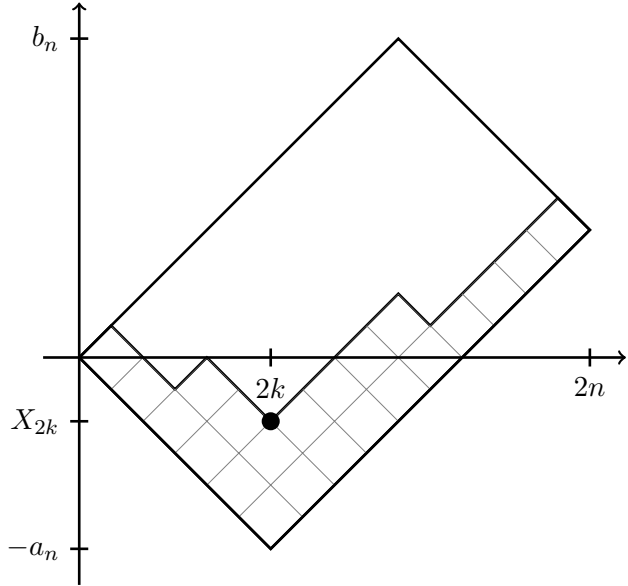

Computing the probability $\pr\paren*{X_k=\ell}$ is straightforward. A lattice
path passing through $(k,\ell)$ is composed of a path from $(0,0)$ to $(k,\ell)$
and a path from $(k,\ell)$ to $(a+b,b-a)$. Since $X_k$ has the same parity as
$k$, the exact formulas also depend on this parity.

\begin{prop}
  \label{prop:1marg}
  The 1-dimensional marginal of $X$ under $\pr$ is given by
  \begin{align}
    \label{eq:8}
    \pr_{a,b}^q(X_{2k}=2i) &= \frac{q^{(k+i)(a -k +i)}}{Z_{a,b}(q)}
    Z_{k-i,k+i}(q) Z_{a-k+i,b-k-i}(q) \intertext{and}
    \label{eq:9}
    \pr_{a,b}^q(X_{2k+1}=2i+1) &= \frac{q^{(k+i+1)(a -k +i)}}{Z_{a,b}(q)}
    Z_{k-i,k+i+1}(q) Z_{a-k +i,b-k-i-1}(q)
  \end{align}
\end{prop}

We focus our attention to the behaviour of the process at even times, which is
sufficient to study its scaling limit, since the process has bounded steps
($+1$/$-1$).


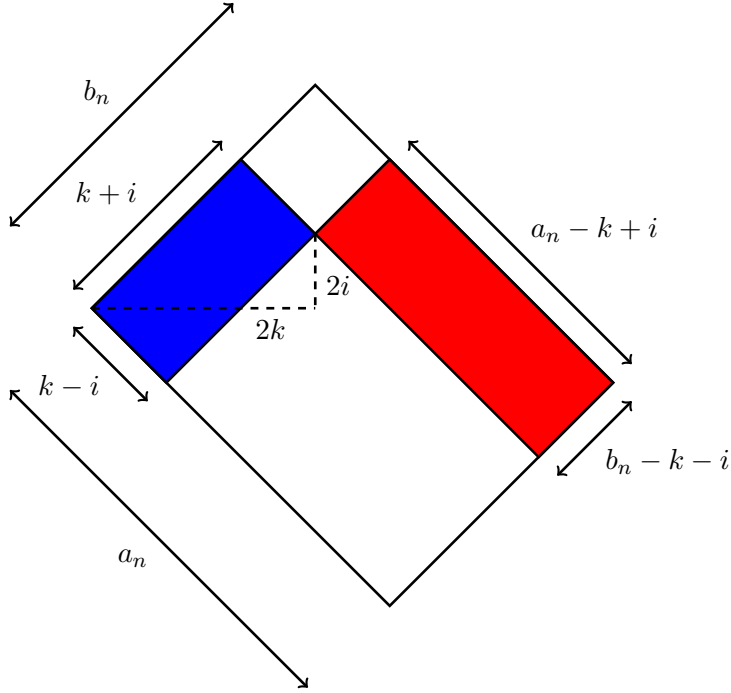
\begin{figure}[htb]
  \centering
  \begin{tikzpicture}[line width=1pt,scale=.7]
    \begin{scope}[rotate=-45,auto]
      \path (0,0) coordinate (origin) (8,6) coordinate (rightcorner) (2,4)
      coordinate (center);
      \draw [fill=blue] (origin) rectangle (center);
      \draw [fill=red] (center) rectangle (rightcorner);
      \draw (origin) rectangle (rightcorner);
      \path (center |- origin) coordinate (p1);
      \draw[<->] (0,-.5) to node [swap] {$k-i$} (0,-.5 -| center);
      \draw[<->] (0,-2.2) to node [swap] {$a_n$} (0,-2.2 -| rightcorner);
      \draw[<->] (-.5,0) to node {$k+i$} (-.5,0 |- center);
      \draw[<->] (-2.2,0) to node {$b_n$} (-2.2,0 |- rightcorner);
      \draw[<->] (0,6.5 -| center) to node {$a_n-k+i$} (0,6.5 -| rightcorner);
      \draw[<->] (8.5,0 |- center) to node [swap] {$b_n-k-i$} (8.5,0 |-
      rightcorner);       
    \end{scope}
    \draw[dashed] (0,0) to node[pos=0.8,below] {$2k$} (0,0 -| center);
    \draw[dashed] (center) to node[pos=0.7,right] {$2i$} (0,0 -| center);
  \end{tikzpicture}
  \caption{Illustration of Equation \eqref{eq:8}.}
  \label{fig:marginalprobability}
\end{figure}

We now mention a unimodality result for the distribution above which turns out
to be useful later.

\begin{lem}
  \label{lem:unimodal}
  The function $\ell \mapsto \pr_{a,b}^q(X_{2k} =2\ell)$ is unimodal: there exists an
  integer $L^{\sharp}(k)=L^{\sharp}_{a,b,q}(k)$ such that
  \begin{equation}
    \frac{\pr_{a,b}^q(X_{2k}=2(\ell+1))}{\pr_{a,b}^q(X_{2k}=2\ell)} \leq 1 \iff \ell \geq
    L^{\sharp}(k). 
    \label{eq:unimodal}
  \end{equation}
\end{lem}
\begin{proof}
  Writing the ratio
  \begin{align*}
    \frac{\pr_{a,b}^q(X_{2k}=2(\ell+1))}{\pr_{a,b}^q(X_{2k}=2\ell)} &=
    \frac{Z_{k-\ell-1,k+\ell+1}(q)
      Z_{a-k+\ell+1,b-k-\ell-1}(q)}{Z_{k-\ell,k+\ell}(q)
      Z_{a-k+\ell,b-k-\ell}(q)}q^{(a-k+\ell+1)(k+\ell+1)-(a-k+\ell)(k+\ell)} \\
    &=\frac{(1-q^{k-\ell})(1-q^{b-k-\ell})}{(1-q^{k+\ell+1)}(1-q^{a-k+\ell+1})}q^{a+2\ell+1}
  \end{align*}
  This ratio is smaller than $1$ if and only if
  \begin{align*}
    (1-q^{k-\ell})(1-q^{b-k-\ell})q^{a+2\ell+1} &\leq (1-q^{k+\ell+1})(1-q^{a-k+\ell+1}) \\
    \iff (q^{a+b+1}-1) &\leq q^{a+1}(q-1)q^{2\ell}
    + \paren*{q^{k+1}(q^a-1)+q^{a-k+1}(q^b-1)}q^\ell
  \end{align*}
  Both terms on the right hand side are increasing continuous functions of $\ell$,
  proving the existence of an integer $L^{\sharp}(k)$ with the asserted
  property.
\end{proof}

\begin{prop}
  \label{prop:2marg}
  The 2-dimensional marginal of $(X_{2k})_{0\leq k\leq n}$ is given by
  \begin{multline*}
    \pr_{a,b}^q(X_{2k}=2i,X_{2l}=2j)=\\
    \frac{Z_{k+i,k-i}(q)Z_{j+l-i-k,l+i-k-j}(q)Z_{a-l+j,b-l-j}(q)}{Z_{a,b}(q)}
    q^{(a-(k-i))(k+i)+(l+j-k-i) (a-l+j)}.
  \end{multline*}
\end{prop}

Now we have exact expressions for probabilistic quantities, we let the dimensions of the box depend on $n$ and investigate the limiting behaviour of the random Young diagram under $\pr_{a_n,b_n}^{q}$ in the regime
\begin{equation*}
  a_n+b_n=2n\to +\infty,\quad \frac{a_n}{2n}\to \rho\in(0,1),\quad -n\log q \to c \in \mathbb{R}.
\end{equation*}

\section{Asymptotics of $\mathbf{q}$-factorials}
In order to study the limit of the process as the size of the box goes to infinity, we first need the asymptotic behavior of the
$q$-factorial.

%
%

\begin{prop}
  \label{prop:qstir2}
  Let $\varepsilon>0$. For $\ell > (1-q)^{-1} \varepsilon$, 
  \begin{equation*}
    (\ell !)_q = \sqrt{2\pi \ell_q} \exp \left( \int_{0}^{\ell} \log t_q \d t \right) \left( 1+O(1-q) \right)
  \end{equation*}

\end{prop}

begin{equation*}
  Consider $\log \ell!_q =
  =\sum_{k=1}^{\ell} \log k_q$.
  By the Euler-Maclaurin formula, we have
  \begin{equation*}
      \log \ell!_q = \int_0^\ell\log(t_q)\d{t} +
      \tfrac{1}{2}\paren*{\log \ell_q}
      + \int_1^\ell
      B_1(t) \frac{\d \log t_q}{\d t}
      -\int_0^1 \log (t_q) \d t.
      \label{eq:2}
  \end{equation*}
  where $B_1(t)=t-\floor{t}-\frac{1}{2}$ is the first periodic Bernoulli polynomial. 

The integral $\int_{0}^{1} \log t_q \d t= 1+O(1-q)$.
For the second integral in \eqref{eq:2}, we add and subtract $\frac{1}{t}$ to $\frac{\d \log t_q}{\d t}$ in order to get
  \begin{equation}
    \label{eq:11}
    \int_1^\ell B_1(t)\frac{\d \log t_q}{\d t}\d{t} 
    = \int_1^\ell B_1(t)\tfrac{1}{t}\d{t} + \int_{1}^{\ell}
    B_1(t)\paren*{\frac{\d \log t_q}{\d t}-\frac{1}{t}}\d{t} 
  \end{equation}
  The function $f(t) = \frac{\d \log t_q }{\d t} - \frac{1}{t}$ is continuous on
  $\mathbb{R}_+$ and its derivative is bounded by an absolute constant $M$ times $(1-q)^2$, i.e. $\abs{f'(t)}\leq M(1-q)^2$ for all $t\in\mathbb{R}_+$. On the
  interval $\brac{k-1,k}$, the function $B_1(t)$ is given by $B_1(t)
  = t-k+\frac{1}{2}$ and has the antiderivative $\tilde{B}_{n,k}(t) = \frac{1}{2}t^2
  -(k-\frac{1}{2})t +\frac{k(k-1)}{2}$ with $\tilde{B}_{n,k}(k-1) =
  \tilde{B}_{n,k}(k) = 0$. By partial integration we get
  \begin{align*}
    \abs*{\int_1^{\ell} B_1(t)f(t)\d{t}} \leq 
    M\sum_{k=1}^{\ell-1} \int_{k}^{k+1} \abs{\tilde{B}_{n,k}(t)}\d{t}
    = \frac{M \ell(1-q)^2}{12}
  \end{align*}
  proving that this term is $O(1-q)$ as $n\to \infty$. The first term of
  \eqref{eq:11} is found to be
  \begin{align}
    \begin{split}
      \int_1^\ell B_1(t)\tfrac{1}{t}\d{t} = \ell-1 -(\ell+\tfrac{1}{2})\log \ell + \log
      \ell! &= -1+\tfrac{1}{2}\log2\pi + O(\tfrac{1}{\ell})
    \end{split}
    \label{eq:stir2}
  \end{align}
  by the classical Stirling approximation $\ell! = \sqrt{2\pi \ell}
  \paren*{\frac{\ell}{e}}^\ell\paren*{1+O(\tfrac{1}{\ell})}$.  Combining \eqref{eq:2},
  \eqref{eq:11} and \eqref{eq:stir2}, we get the result. 

Using the notations of Section \ref{sec:present}, we restate the q-Stirling's formula in a convenient form for this particular context.

\begin{cor}[q-Stirling's Formula]
  \label{prop:qstirling}
  Let $c\in\mathbb{R_+^*}$ and fix $\epsilon>0$
  In the limit when $n$ goes to infinity, with $q=e^{-\frac{c}{n}}$, the following
  asymptotics hold for all $\ell > \varepsilon n$
  \begin{equation*}
    \ell!_q = \sqrt{2\pi n}\sqrt{\frac{e^{\frac{c\ell}{n}}-1}{c}} n^\ell
    \exp(n S_c(\tfrac{\ell}{n})) \paren[\big]{1+O(\tfrac{1}{\ell})}, 
  \end{equation*}
  where
  \begin{equation*}
    S_c(\alpha)= \int_{0}^{\alpha} \log\paren*{ \frac{1-e^{-cx}}{c} } \d{x}
    = \alpha S_{\alpha c}(1).
  \end{equation*}
  
  In particular, there exist positive constants $m$ and $M$ such that for all $n$, and all
  $\ell$ between 1 and $n$,
  \begin{equation*}
    m \sqrt{2\pi n}
    \sqrt{\frac{e^{\frac{c\ell}{n}}-1}{c}}    {n}^\ell 
    \exp(n S_c(\tfrac{\ell}{n})) < \ell!_q < M \sqrt{2\pi
      n}\sqrt{\frac{e^{\frac{c\ell}{n}}-1}{c}} {n}^\ell 
    \exp(n S_c(\tfrac{\ell}{n})) .
  \end{equation*}
\end{cor}


We turn now to the the asymptotics of the 1- and 2-dimensional marginals of $X$ under
$\pr$ given by Propositions \ref{prop:1marg} and \ref{prop:2marg}. For this purpose, define
\begin{equation*}
  f_c(x,y) = S_c(x+y) -S_c(x)-S_c(y) \quad\text{ and }\quad
  h_c(x,y) = \sqrt{\frac{c(e^{c(x+y)}-1)}{(e^{cx}-1)(e^{cy}-1)}}.
\end{equation*}
If $x=\frac{\ell}{n}$ and $y=\frac{k}{n}$, the asymptotics of the $q$-binomial coefficient
is given by
\begin{equation*}
  Z_{l,k}^{(q)}=\binom{\ell+k}{\ell}_q = \sqrt{\frac{1}{2\pi n}}
  h_c(x,y)\exp(nf_c(x,y))\paren[\big]{1+O(\tfrac{1}{n})}.
\end{equation*}

As a consequence,
\begin{cor} \label{cor:1-marg} For all $\varepsilon>0$, for all $k$ between $\varepsilon
  n$ and $(1-\varepsilon)n$, with $s=k/n$ and $x=i/n$, we have
  \begin{equation*}
    \pr(X_{2k}=2i) = \sqrt{\frac{1}{2\pi n}} H^{(1)}_{\rho,c}(s,x)
    \exp(n F^{(1)}_{\frac{a_n}{2n},c}(s,x))\paren[\big]{1+O(\tfrac{1}{n})},
  \end{equation*}
  and there exists a constant $M>0$ such that for all $n$, $i$ and $k$,
  \begin{equation*}
    \pr(X_{2k}=2i) \leq M \sqrt{\frac{1}{2\pi n}}
    H^{(1)}_{\rho,c}(s,x) \exp(n
    F^{(1)}_{\rho,c}(s,x)),
  \end{equation*}
  where
  \begin{align*}
    F^{(1)}_{\rho,c}(s,x) &= -c(2\rho-s+x)(s+x) + f_c(s-x,s+x) + f_c(2\rho-s+x,2-2\rho-s-x) - f_c(2\rho,2-2\rho),\\
    H^{(1)}_{\rho,c}(s,x)&=\frac{h_c(s+x,s-x)h_c(2\rho-s+x,2-2\rho-s-x)}{h_c(2\rho,2-2\rho)}.
  \end{align*}

\end{cor}
In the first equality, when replacing $\frac{a_n}{2n}$ by $\rho$ in the indices of
$H^{(1)}$, we make an error of order $O(\frac{1}{n})$, which is absorbed in the factor
$(1+O(\frac{1}{n}))$.  Making the same substitution in the indices of $F^{(1)}$ would
change the multiplicative constant in the asymptotics.

\begin{cor}\label{cor:2-marg}
  For all $\varepsilon>0$, for all $k < l$ between $\varepsilon n$ and
  $({1}-\varepsilon)n$ then, with $s=k/n$, $t=l/n$ and $x=i/n$, $y=j/n$
  \begin{equation*}
    \pr(X_{2k}=2i,X_{2l}=2j) = \frac{1}{2\pi
      n} H^{(2)}_{\rho,c}(s,t,x,y)
    \exp(n F^{(2)}_{\frac{a_n}{2n},c}(s,t,x,y))(1+O(\tfrac{1}{n})), 
  \end{equation*}
  where
  \begin{align*}
    F^{(2)}_{\rho,c}(s,t,x,y) &= -c\bigl((2\rho-s+x)(s+x)+(t-s+y-x)(2\rho-t+y)\bigr) \\&
    \qquad+ f_c(s-x,s+x) + f_c(t-s+y-x,t-s-y+x)\\
    &\qquad +f_c(2\rho-t+y,2-2\rho-t-y) - f_c(2\rho,2-2\rho),\\
    H^{(2)}_{\rho,c}(s,t,x,y)
    &=\frac{h_c(s+x,s-x)h_c(t-s+y-x,t-s-y+x)h_c(2\rho-t+y,2-2\rho-t-y)}{h_c(2\rho,2-2\rho)}.
  \end{align*}
\end{cor}

We give now a slight refinement of Corollary~\ref{cor:2-marg} that will be useful in the
fluctuations.

\begin{cor}
  For all $\varepsilon>0$, for all $k < l$ between $\varepsilon n$ and $(1-\varepsilon)n$
  then, with $s=k/n$, $t=l/n$ and $x=i/n$, $y=j/n$, for all $s'=s+o(1)$, $t'=t+o(1)$,
  $x'=x+o(1)$, $y'=y+o(1)$:
  \begin{equation*}
    \pr_n(X_{2k}=2i,X_{2l}=2j) = \frac{1}{2\pi
      n} H^{(2)}_{\rho,c}(s',t',x',y')
    \exp(n F^{(2)}_{\frac{a_n}{2n},c}(s,t,x,y))(1+o(1)).
  \end{equation*}
\end{cor}

\section{Limit Shape}

We associate to the lattice path $(X_k)$ the continuous piecewise linear
function, $X:s\mapsto X_s$ defined on $[0,2n]$, which coincides with $X_k$ when $s=k$.
The graph of the function $X$ is the boundary of the random Young diagram we consider.


Let $L_{\rho,c}$ be the function on $[0,{1}]$ defined by
\begin{align}
  \forall t\in[0,{1}],\quad L_{\rho,c}(t) &=
  1-2\rho+\frac{1}{c}\log
  \frac{\sinh(ct)+e^{c(2\rho-1)}\sinh(c({1}-t))} {\sinh {c}\nonumber}
  \\
  \label{eq:7}
   &= \tfrac{1}{c} \log\frac{e^{-ct} -e^{ct} +e^{-c(2-2\rho-t)}
    -e^{-c(t-2\rho)}} {e^{-c(2-2\rho)}-e^{2c\rho}}. 
\end{align}

This is the limit shape of the rescaled random process, in the sense made
precise in Theorem \ref{thm:1} below. First we need to estimate the proximity between the value of the function $L_{\rho,c}$ and the most probable value for $X_{2k}$ denoted by $L^{\sharp}_n(k)=L_{a_n,b_n,e^{-c/n}}(k)$ in Lemma \ref{lem:unimodal}.

\begin{lem}
  \label{lem:limitshape}
  For all $n$ and all $0\leq k\leq n$,
  $ \abs*{\tfrac{1}{n}L^{\sharp}_{n}(k)-L_{\rho,c}(\tfrac{k}{n})} \leq \tfrac{1}{n}
  + \abs*{\tfrac{a_n}{2n}-\rho}$, where $q=e^{-\frac{c}{n}}$.
 
  Consequently, in the limit when $n$ goes to infinity, and $a_n/2n$ goes to $\rho$,
  \begin{equation}\label{eq:41}
    \abs*{\tfrac{1}{n}L^{\sharp}_{n}(k)-L_{\rho,c}(\tfrac{k}{n})} \rightarrow 0.
  \end{equation}
\end{lem}

\begin{proof}
  The ratio of probabilities \eqref{eq:unimodal} used to define $L^{\sharp}(k)=L^{\sharp}_{n}(k)$ can
  be rewritten in terms of the function
  \begin{equation*}
   (\rho,c,t,x)\mapsto R_{\rho,c}(t,x) = e^{-c} \frac{\sinh(\frac{c}{2}(t-x))
      \sinh(\frac{c}{2}(2-2\rho-t-x))} {\sinh(\frac{c}{2}(t+x))
      \sinh(\frac{c}{2}(2\rho-t+x))},
  \end{equation*}
  namely as $R_{\frac{a_n}{2n},c}(\tfrac{k}{n}, \frac{\ell}{n})$. The function
  $x\mapsto R_{\rho,c}(t,x)$ is decreasing, and by definition $L^{\sharp}_n(k)$
  is the smallest integer $\ell$ such that
  $R_{\tfrac{a_n}{2n},c}(\tfrac{k}{n},\tfrac{\ell}{n}) \leq 1$.

  On the other hand, a computation shows that for each $t\in
  \brac*{0,{1}}$, the equation
  \begin{equation*}
    R_{\rho,c}(t,x) = 1
  \end{equation*}
  has the (unique) solution $x=L_{\rho,c}(t)$. Substituting $\tfrac{a_n}{2n}$ for
  $\rho$, we conclude that
  \begin{equation}
    \label{eq:43}
    \abs*{\tfrac{1}{n}L^{\sharp}_{n}(k)-L_{\tfrac{a_n}{2n},c}(\tfrac{k}{n})} \leq
    \tfrac{1}{n} 
  \end{equation}
  for all $n$ and all $k\leq n$.

  Differentiating $L_{\rho,c}$ with respect to $\rho$, we find that
  \begin{equation*}
    \abs*{\dd{L_{\rho,c}}{\rho}(t)} \leq 1
  \end{equation*}
  for all $t$, and so by the mean value theorem,
  \begin{equation}
    \label{eq:44}
    \abs*{L_{\tfrac{a_n}{2n},c}(\tfrac{k}{n}) - L_{\rho,c}(\tfrac{k}{n})} \leq
    \abs*{\frac{a_n}{2n} -\rho}
  \end{equation}
  for all $k\leq n$. From \eqref{eq:43}, \eqref{eq:44} and the
  triangle inequality, we get the expected result.
\end{proof}

\begin{rem}
  It turns out that $\log R_{\rho,c}(t,x)$ coincides with the derivative of $F^{(1)}_{\rho,c}(s,x)$ with respect to $x$. As a consequence, 
  $L_{\rho,c}(s)$ can also be viewed as the argmax of $F^{(1)}_{\rho,c}(s,\cdot)$, which is nonpositive and vanishes at $L_{\rho,c}(s)$. See Corollary \ref{cor:F(s,L(s))=0}.
\end{rem}


\begin{thm}
  \label{thm:1}
  The boundary of the rescaled random Young diagram converges in probability, for the uniform
  topology, to the curve of $t\mapsto L_{\rho,c}(t)$.
  \begin{equation*}
    \forall \varepsilon>0,\ \lim_{n\rightarrow 0}
    \mathbb{P}^{e^{-c/n}}_{a_n,b_n}\paren[\Bigg]{\sup_{t\in\brac*{0,1}}
      \abs*{\tfrac{1}{2n}X_{2tn}-L_{\rho,c}(t)} > \varepsilon} = 0, 
  \end{equation*}
where
\begin{equation*}
 L_{\rho,c}(t) = \tfrac{1}{c} \log\frac{e^{-ct} -e^{ct} +e^{-c(2-2\rho-t)}
    -e^{-c(t-2\rho)}} {e^{-c(2-2\rho)}-e^{2c\rho}}.
\end{equation*}

\end{thm}

Remark that in the case of a square box ($\rho=1/2$), the expression for $L_{\rho,c}$ boils down to:
\begin{equation*}
 L_{\frac{1}{2},c}(t)=  \frac{1}{c} \log \frac{\cosh(c(t-\frac{1}{2}))}{\cosh \frac c 2}
\end{equation*}

\begin{proof}
    Fix $\varepsilon>0$. For $t<\frac{\varepsilon}{2}$ or
  $t>1-\frac{\varepsilon}{2}$, the difference
  $\abs*{\frac{1}{2n}X_{2tn}- L{\rho,c}(t)}$ is always smaller than $\varepsilon$. We
  have to control what happens for $t\in(\frac{\varepsilon}{2},1-
  \frac{\varepsilon}{2})$. Using the fact that $L_{\rho,c}$ is differentiable, and its
  derivative with respect to $t$ is bounded by $1$, we have that
  \begin{equation*}
    \abs*{L_{\rho,c}\paren*{\tfrac{1}{n}\lfloor tn \rfloor} - L_{\rho,c}(t)} \leq \frac{1}{n},
  \end{equation*}
  which is smaller than $\frac{\varepsilon}{3}$ for $n$ sufficiently large. The
  same is true of $\abs*{\frac{1}{2n}X_{2tn}-\frac{1}{2n}X_{2\lfloor
      tn\rfloor}}$. Thus, by an $\frac{\epsilon}{3}$-argument, to control the
  sup over $(\frac{\varepsilon}{2},1 -\frac{\varepsilon}{2})$, it is
  sufficient to control what happens at points of the form $t=\frac{k}{n}$:
  \begin{multline}
    \pr\paren*{\sup_{t\in(\frac{\varepsilon}{2}, 1-\frac{\varepsilon}{2})}
    \abs*{\tfrac{1}{2n}X_{2tn}-L_{\rho,c}(t)}>\varepsilon} \\ \leq
    \sum_{k\in\mathbb{N}\cap n(\frac{\varepsilon}{2}, 1 -
    \frac{\varepsilon}{2})} \pr\paren[\big]{X_{2k} > 2n (L_{\rho,c}(\tfrac{k}{n})+
    \varepsilon) \ \text{or}\ X_{2k} < 2n(L_{\rho,c}(\tfrac{k}{n})-\varepsilon)}
    \label{eq:bound_limit_shape}
  \end{multline}
  But it follows from \eqref{eq:41}, that for $n$ sufficiently large,
  $2n\paren[\big]{L(\frac{k}{n})+\varepsilon} \geq 2L^{\sharp}_{n}(k)$ for all
  $k\in\mathbb{N}\cap
  n(\frac{\varepsilon}{2},1-\frac{\varepsilon}{2})$, and thus using
  the unimodality of the law of $X_{2k}$
  \begin{align*}
    \pr\paren[\big]{X_{2k} > 2n(L_{\rho,c}(\tfrac{k}{n})+\varepsilon)}
    &= \sum_{l>n\paren{L_{\rho,c}(\frac{k}{n})+\varepsilon}} \pr(X_{2k} =2l) \\
    &\leq n\times \pr\paren*{X_{2k}=2\floor*{n(L_{\rho,c}(\tfrac{k}{n})+\varepsilon)}}.
  \end{align*}
  which by Corollary \ref{cor:1-marg} is exponentially small, uniformly in $k$,
  as $n$ goes to infinity. Thus for $n$ sufficiently large, the sum on the RHS
  of \eqref{eq:bound_limit_shape} can be made smaller than any positive number.
\end{proof}

\begin{rem}
  \label{rem:1}
  From the proof and Lemma \ref{lem:limitshape}, we see that the convergence
  statement in the theorem holds uniformly for a family of sequences $(a_n)$ and
  $(b_n)$, so long as the convergence $\frac{a_n}{2n}\to \rho$ is uniform for the
  family. 
\end{rem}

\section{Fluctuations}

We now study the fluctuations of the interface around the limit shape. We define
for $n\in\mathbb{N}^*$ a new rescaled process:
\begin{equation}
  \forall t\in\brac*{0,1},\quad \tilde{X}_t = \tilde{X}^{(n)}_t =
  \sqrt{n}\paren*{\frac{1}{2n}X_{2nt}- L_{\rho,c}(t)}. 
  \label{eq:proc_fluct}
\end{equation}

We place ourselves in the space $D$ of c\`ad-l\`ag paths on $[0,1]$ endowed
with its usual topology. We state now our second main result for the
convergence of the fluctuations of the interface to the Ornstein Uhlenbeck
bridge. See the appendix for the definition and some properties of this process.

Define 
\begin{equation}
 f(s)=\frac{1}{\sqrt{2}} \frac{\sqrt{(e^{2c\rho}-1)(1-e^{-2c(1-\rho)})}} {\sinh
      cs + e^{c(2\rho-1)}\sinh(c(1 -s))}
      =\frac{\sqrt{2\sinh(c\rho)\sinh(c(1-\rho))}}{e^{c(\frac{1}{2}-\rho)}\sinh(cs)+e^{c(\rho-\frac{1}{2})}\sinh(c(1-s))}.
      \label{eq:def_f}
\end{equation}

\begin{thm}\label{thm:2}
  The sequence $(\tilde{X}^{(n)}_s / f(s))_n$, converges weakly in $D$ to the
  Ornstein-Uhlenbeck bridge $(Y_s)_{s\in[0,1]}$, which is the Gaussian process on $[0,1]$ 
  with covariance 
\begin{equation*}
 \E\brac{Y_s Y_t} =  \frac{\sinh(cs)\sinh(c(1-t))}{c\sinh(c)},
\end{equation*}
for $0\leq s<t\leq {1}$.

\end{thm}

\subsection{Two-point correlations}
\label{subsec:fluc2pt}

To prove convergence of the two-dimensional marginal to a Gaussian process, we
apply a saddle-point method, ie. we need to show that the function $F^{(2)}$,
which governs the exponential decay of 2-dimensional marginals (see Corollary
\ref{cor:2-marg}), has a critical point `on the limit shape', and that it takes
its maximal value of $0$ at this point.  We will prove that this is indeed the
case, but first we need a lemma describing the limit shape in a subrectangle.

\begin{lem}
  \label{lem:1}
  The limit shape $L$ satisfies the relations
  \begin{align}
    \label{eq:12}
    \frac{1}{1-s}\paren[\big]{L_{\rho,c}(t)-L_{\rho,c}(s)} &=
    L_{\rho',c'}\paren[\big]{\tfrac{t-s}{1-s}} \shortintertext{and}
    \label{eq:27}
    \frac{1}{t}L_{\rho,c}(s) &= L_{\rho'',c''}\paren*{\tfrac{s}{t}}
  \end{align}
  for $0<s<t<{1}$, with
  \begin{alignat*}{2}
    \rho' &= \frac{2\rho-s+L_{\rho,c}(s)}{2(1-s)}, \qquad\qquad& c'&=c(1-s)
    \shortintertext{and} \rho'' &= \frac{t-L_{\rho,c}(t)}{2t}, & c'' &= tc.
  \end{alignat*}
\end{lem}
\begin{proof}
  These relations can be checked analytically. We sketch now a less computational argument for Equation \eqref{eq:27}.
The same strategy applies for Equation \eqref{eq:12}.
Fix $t\in(0,1)$, and take a sequence of boxes with sidelengths $[0,a_n]\times[0,b_n]$ (with $a_n+b_n=2n$, $a_n/2n\rightarrow \rho$).
Theorem \ref{thm:1} states that, when $n$ goes to infinity, the probability under $\pr_{a_n,b_n}^{e^{-c/n}}$ that
 $X_{nt}/n$ converges to $L_{\rho,c}(t)$ goes to 1. As a consequence, the restriction of the limit shape $L_{\rho,c}$ to the,
 rescaled by $n$, lower left subbox $[0,\frac{t-L_{\rho,c}(t)}{2}]\times[0,\frac{t+L_{\rho,c}(t)}{2}]$ is the rescaled limit shape for boxes with ratio
  $\frac{t-L_{\rho,c}(t)}{t-L_{\rho,c}(t)+t+L_{\rho,c}(t)}=\frac{t-L_{\rho,c}(t)}{2t}$. The rescaling factor is the perimeter of
the subbox $t$. The value of the parameter $q$ remains the same. Since $q=e^{-c/n} = e^{-\frac{c''}{nt}}$, we get $c''=tc$.
\end{proof}

This lemma basically describes the limit shape of the process restricted to a
subrectangle defined by a point $(s,L_{\rho,c}(s))$ and either the right or left
corner of the original rectangle.

\begin{lem}
  For all $0<s<t<1$, the functions $F^{(1)}$ and $F^{(2)}$ satisfy the relations
  \begin{align}
    \label{eq:30}
    F^{(2)}_{\rho,c}(s,t,x,y) &= F^{(1)}_{\rho,c}(s,x) + (1-s)
    F^{(1)}_{\rho',c'}(t',y') \shortintertext{and}
    \label{eq:17}
    F^{(2)}_{\rho,c}(s,t,x,y) &=
    F^{(1)}_{\rho,c}(t,y)+ t F^{(1)}_{\rho'',c''}(s'',x'')
  \end{align}
  with
  \begin{alignat*}{4}
    \rho' &= \frac{2\rho-s+x}{2(1-s)},\qquad & t' &=\frac{t-s}{1-s}, \qquad &
    y'&=\frac{y-x}{1-s}, \qquad & c'&=c(1-s) \shortintertext{and} \rho''
    &=\frac{t-y}{2t},\qquad & s''&=\frac{s}{t},\qquad & x''
    &=\frac{x}{t},\qquad & c''&=tc.
  \end{alignat*}
\end{lem}
\begin{proof}
  The relation \eqref{eq:30} follows from the identity 
  \begin{equation*}
    \pr(X_{2k}=2i,X_{2\ell}=2j) = \pr (X_{2k}=2i) \pr(X_{2\ell}=2j| X_{2k}=2i).
  \end{equation*}
  allied to scaling arguments similar to those used in the proof of the previous lemma.
Equation \eqref{eq:17} follows from the conditioning on $X_{2l}$
  instead of $X_{2k}$.
\end{proof}

\begin{lem}
  \label{lem:ddF1}
  The partial derivatives of $F^{(1)}$ with respect to $\rho$,
   $s$ and $x$ all
  vanish at $x=L_{\rho,c}(s)$.
\end{lem}

\begin{cor}
  \label{cor:ddF2}
  The partial derivatives of $F^{(2)}$ with respect to $\rho$, $s$, $t$, $x$ and
  $y$ all vanish at $x=L_{\rho,c}(s)$ and $y=L_{\rho,c}(t)$.
\end{cor}
\begin{proof}
  For the partial derivatives with respect to $\rho$, use \eqref{eq:17}. For
  $\dd{F^{(2)}}{s}$ and $\dd{F^{(2)}}{x}$, use \eqref{eq:17} together with
  \eqref{eq:27}. Similarly, for $\dd{F^{(2)}}{t}$ and $\dd{F^{(2)}}{y}$, use
  \eqref{eq:30} together with \eqref{eq:12}.
\end{proof}

\begin{cor}
  \label{cor:F(s,L(s))=0}
  $F^{(1)}_{\rho,c}(s,L_{\rho,c}(s))=0$ and
  $F^{(2)}_{\rho,c}(s,t,L_{\rho,c}(s),L_{\rho,c}(t))=0$.
\end{cor}
\begin{proof}
  By the chain rule and Lemma \ref{lem:ddF1} we get
  \begin{equation*}
    \frac{d}{ds}(F^{(1)}_{\rho,c}(s,L_{\rho,c}(s))) = 0.
  \end{equation*}
  Taking the limit $s\to 0$ yields zero, proving the claim. For the second
  identity, use the first one together with \eqref{eq:17} and \eqref{eq:27}.
\end{proof}

\begin{prop} \label{prop:fluct2dim}
  Assume $\frac{a_n}{2n}=\rho+o\paren[\big]{\frac{1}{\sqrt{n}}}$. Let $s<t$ in
  $\brac*{0,1}$. The joint law of
  $(\tilde{X}^{(n)}_s,\tilde{X}^{(n)}_t)$ converges to the law of the
  centered
  2-dimensional Gaussian vector
  \begin{equation*}
    (\tilde{X}_s,\tilde{X}_t)=(f(s) Y_s, f(t) Y_t),
  \end{equation*}
  where $f(s)$ is defined in as in Equation \eqref{eq:def_f},
  and $(Y_t)_{t\in[0,1]}$ is an Ornstein-Uhlenbeck bridge on the interval
  $[0,1]$ with parameter $c$ (see Appendix).
  
\end{prop}

\begin{rem}
  In the case of a square ($\rho=\frac{1}{2}$), the expression of $f$ simplifies
  drastically to
  \begin{equation*}
    f(x)= \frac{1}{\sqrt{2}\cosh\big( c\bigl(s-\frac{1}{2}\bigr)\bigr)}.
  \end{equation*}
\end{rem}

\begin{proof}
  Let $a_1, b_1, a_2, b_2$ be four real numbers such that $a_1<b_1$ and
  $a_2<b_2$. In the following, we write $L$ for $L_{\rho,c}$. Clearly,
  $\tilde{X}^{(n)}_s \in \brac{a_1,b_1}$ if and only if $X_{2ns}$ is in the
  interval
  \begin{equation*}
    2nL(s) +\brac*{2a_1\sqrt{n},2b_1\sqrt{n}}.
  \end{equation*}
  Since $\abs*{X_{2ns}-X_{2\floor{ns}}} \leq 2$ and $\abs*{X_{2nt}-X_{2\floor{nt}}} \leq 2$,
we just need to compute the limit of
  \begin{equation*}
    \pr\paren*{\paren*{X_{2\floor{ns}}, X_{2\floor{nt}}} \in A_1\times A_2}.
  \end{equation*}
  as $n$ goes to infinity.

  Set $k_n =\floor{ns}$ and $\ell_n=\floor{nt}$. One has:
  \begin{multline*}
    \pr\paren*{\paren*{X_{2\floor{ns}}, X_{2\floor{nt}}} \in A_1\times A_2}=\sum_{(i,j)\in \frac{1}{2}A_1\times \frac{1}{2}A_2}
    \pr\paren*{X_{2k_n} =
      2i, X_{2\ell_n} = 2j}\\
    = \sum_{(i,j)\in \frac{1}{2}A_1\times \frac{1}{2}A_2} \frac{1}{2\pi n}
    H^{(2)}_{\rho,c}(s,t,L(s),L(t)) \exp\paren*{n F^{(2)}_{\frac{a_n}{2n},c}
      \paren*{\tfrac{k_n}{n},\tfrac{\ell_n}{n},
        \tfrac{i}{n},\tfrac{j}{n}}}\paren[\big]{1+o(1)}.
  \end{multline*}

    
  Here, $F^{(2)}\leq 0$ everywhere and the number of terms is $O(n)$, so the
  error terms $o(1)$ can be replaced with an $o(1)$ term outside the sum. Thus,
  we are left with a Riemann sum for the double integral
  \begin{equation}
    \label{eq:21}
    \frac{1}{2\pi n} H^{(2)}_{\rho,c}(s,t,L(s),L(t))
    \int_{nL(s)+a\sqrt{n}}^{nL(s)+b\sqrt{n}}
    \int_{nL(t)+c\sqrt{n}}^{nL(t)+d\sqrt{n}}
    \exp\brac*{nF^{(2)}_{\frac{a_n}{2n},c}\paren*{\tfrac{k_n}{n},
        \tfrac{\ell_n}{n},\tfrac{u}{n},\tfrac{v}{n}}} \d{v}\d{u}  
  \end{equation}
  in which we make the substitution
  \begin{equation*}
    u=x\sqrt{n}+nL(s),\qquad v=y\sqrt{n}+nL(t)
  \end{equation*}
  to get
  \begin{equation*}
    \frac{1}{2\pi}H^{(2)}_{\rho,c}(s,t,L(s),L(t))
    \int_a^b\int_c^d \exp\brac*{nF^{(2)}_{\frac{a_n}{2n},c}\paren*{\tfrac{k_n}{n},
        \tfrac{\ell_n}{n},\tfrac{x}{\sqrt{n}}+L(s),
        \tfrac{y}{\sqrt{n}}+L(t)}} \d{y}\d{x} 
  \end{equation*}
  We make a second degree Taylor expansion of the function
  $F^{(2)}_{\frac{a_n}{2n},c}\paren*{\frac{k_n}{n},\frac{\ell_n}{n},\cdot,\cdot}$
  at the point $(L(s),L(t))$. Let $z_n$ be the point
  $(\frac{k_n}{n},\frac{\ell_n}{n},L(s),L(t))$ in $\mathbb R^4$.
  \begin{align}
    \label{eq:24}
    \begin{split}
      &nF^{(2)}_{\frac{a_n}{2n},c} \paren*{\tfrac{k_n}{n},\tfrac{\ell_n}{n},
        \tfrac{x}{\sqrt{n}}+L(s), \tfrac{y}{\sqrt{n}}+L(t)} \\
      &\qquad = nF^{(2)}_{\frac{a_n}{2n},c}(z_n) +
      \dd{F^{(2)}_{\frac{a_n}{2n},c}}{x}(z_n)x\sqrt{n} +
      \dd{F^{(2)}_{\frac{a_n}{2n},c}}{y}(z_n)y\sqrt{n} \\ &\qquad\qquad +
      \frac{1}{2}\ddtwo{F^{(2)}_{\frac{a_n}{2n},c}}{x}(z_n)x^2 +
      \frac{1}{2}\ddtwo{F^{(2)}_{\frac{a_n}{2n},c}}{y}(z_n)y^2 +
      \ddd{F^{(2)}_{\frac{a_n}{2n},c}}{x}{y}(z_n)xy +
      O\paren[\big]{n^{-\frac{1}{2}}}
    \end{split}
  \end{align}
  
  We must find the limit of this expression as $n\to\infty$. The first term on
  the right has limit zero, which can be seen by Taylor expansion around
  $(s,t,L(s),L(t))$, using Corollaries \ref{cor:ddF2} and \ref{cor:F(s,L(s))=0}.
  
  Since $F^{(2)}$ is an analytic function, we have in general
  \begin{equation}\label{eq:25}
    F^{(2)}_{\rho+o\paren[\big]{n^{-\frac{1}{2}}},c} \paren*{s+o\paren[\big]{n^{-\frac{1}{2}}},
      t+o\paren[\big]{n^{-\frac{1}{2}}} ,x,y} = F^{(2)}_{\rho,c}(s,t,x,y)
    +o\paren[\big]{n^{-\frac{1}{2}}}
  \end{equation}
  and similarly for all partial derivatives of $F^{(2)}$. Since $\frac{a_n}{2n} =
  \rho+o\paren[\big]{n^{-\frac{1}{2}}}$, $\frac{k_n}{n} = s+o\paren[\big]{n^{-\frac{1}{2}}}$ and
  $\frac{\ell_n}{n} = t+o\paren[\big]{n^{-\frac{1}{2}}}$, this shows that the second and
  third terms on the right in \eqref{eq:24} tend to zero as $n\to\infty$. Thus,
  if $z=(s,t,L(s),L(t))$,
  \begin{multline}
    \label{eq:26}
    nF^{(2)}_{\frac{a_n}{2n},c} \paren*{\tfrac{k_n}{n},\tfrac{\ell_n}{n},
      \tfrac{x}{\sqrt{n}}+L(s), \tfrac{y}{\sqrt{n}}+L(t)}\to
    \frac{1}{2}\ddtwo{F^{(2)}_{\rho,c}}{x} (z) x^2 +
    \frac{1}{2}\ddtwo{F^{(2)}_{\rho,c}}{y} (z) y^2 +
    \ddd{F^{(2)}_{\rho,c}}{x}{y} (z) xy
  \end{multline}
  as $n\to\infty$. From \eqref{eq:25} it also follows that the sequence
  $F^{(2)}_{\frac{a_n}{2n},c}(\frac{k_n}{n},\frac{\ell_n}{n},\cdot,\cdot)$ is
  equicontinuous, hence by Ascoli's theorem that the sequence converges
  uniformly, hence that the double integral \eqref{eq:21} converges to
  \begin{equation*}
    \frac{1}{2\pi} H^{(2)}_{\rho,c}(s,t,L(s),L(t))\int_a^b\int_c^d
    \exp\brac*{\frac{1}{2}\ddtwo{F_{\rho,c}^{(2)}}{x}(z) x^2 +
      \frac{1}{2}\ddtwo{F^{(2)}_{\rho,c}}{y}(z) y^2 +
      \ddd{F^{(2)}_{\rho,c}}{x}{y} (z) xy} \d{y}\d{x}
  \end{equation*}
  To find these double derivatives, we use \eqref{eq:30} and
  \eqref{eq:17}. Hence we need
  \begin{equation*}
    \ddtwo{F^{(1)}_{\rho,c}}{x}(s,L(s)) = c(1-e^{-2c}) e^{2cs}
    \frac{(1-e^{-c(2-2\rho)}-e^{-2cs}+e^{-c(2s-2\rho)})^2}
    {(1-e^{2c\rho})(1-e^{-2cs})(1-e^{-c(2-2\rho)})(1-e^{-c(2-2s)})} 
  \end{equation*}
  and exploit the fact that
  \begin{align*}
    \ddtwo{F^{(2)}_{\rho,c}}{x}(s,t,L(s),L(t)) &=
    \frac{1}{t}\ddtwo{F^{(1)}_{\rho'',c''}}{x}(s'',L_{\rho'',c''}(s'')) \\
    \ddtwo{F^{(2)}_{\rho,c}}{y}(s,t,L(s),L(t)) &=
    \frac{1}{1-s}\ddtwo{F^{(1)}_{\rho',c'}}{x}(t',L_{\rho',c'}(t'))
  \end{align*}
  This way, we find the double derivatives of
  $F^{(2)}_{\rho,c}(s,t,\cdot,\cdot)$, evaluated at the critical point
  $(x,y)=(L(s),L(t))$, to be
  \begin{align*}
    \ddtwo{F^{(2)}_{\rho,c}}{x}(s,t,L(s),L(t)) &=
    \frac{c(1-e^{2ct})(1-e^{2cs}+e^{-c(2-2\rho-2s)}-e^{2c\rho})^2}
    {(1-e^{2cs})(e^{2ct}-e^{2cs})(1-e^{2c\rho})(1-e^{-c(2-2\rho)})}\\
    &= \frac{-c \sinh ct}{\sinh cs \sinh c(t-s) f(s)^2},\\
    \ddtwo{F^{(2)}_{\rho,c}}{y}(s,t,L(s),L(t)) &=
    \frac{c(e^{2cs}-e^c)(1-e^{2ct}+e^{-c(2-2\rho-2t)}-e^{2c\rho})^2}
    {(e^{2cs}-e^{2ct})(e^c-e^{2ct})(1-e^{2c\rho})(1-e^{-c(2-2\rho)})}\\
    &= \frac{-c \sinh c(1-s)}{\sinh c(1-t) \sinh c(t-s) f(t)^2}\\
    \ddd{F^{(2)}_{\rho,c}}{x}{y}(s,t,L(s),L(t)) &=
    \frac{c(1-e^{2cs}+e^{-c(2-2\rho-2s)}-e^{2c\rho})
      (1-e^{2ct}+e^{-c(2-2\rho-2t)}-e^{2c\rho})}
    {(e^{2cs}-e^{2ct}) (1-e^{2c\rho}) (1-e^{-c(2-2\rho)})} \\
    &= \frac{c}{\sinh c(t-s) f(s) f(t)}
  \end{align*}
  The last one, the mixed derivative, was calculated from scratch, i.e. by
  calculating
  \begin{equation*}
    \ddd{F^{(2)}_{\rho,c}}{x}{y} (s,t,x,y)= \frac{c(e^{2ct}-e^{2cs})}
    {(e^{c(t+y)}-e^{c(s+x)})(e^{c(t-y)}-e^{c(s-x)})} 
  \end{equation*}
  and evaluating at the critical point $(x,y)=(L(s),L(t))$. Thus, the Hessian of
  $F^{(2)}_{\rho,c}(s,t,\cdot,\cdot)$ at the critical point is
  \begin{multline*}
    H(F^{(2)}_{\rho,c})= \frac{-c}{\sinh{c s} \sinh{c(1-t)}
      \sinh{c(t-s)}}\times\\ \begin{pmatrix} \frac{1}{f(s)} & 0 \\ 0 &
      \frac{1}{f(t)} \end{pmatrix}
    \begin{pmatrix}
      \sinh ct \sinh c(1-t) & -\sinh cs \sinh c(1-t)\\
      - \sinh cs \sinh c(1-t) & \sinh c s \sinh c(1-s)
    \end{pmatrix}
    \begin{pmatrix} \frac{1}{f(s)} & 0 \\ 0 & \frac{1}{f(t)} \end{pmatrix}.
  \end{multline*}
  The covariance matrix of the limiting Gaussian distribution is the negative of
  the inverse of $H(F^{(2)})$, which we compute to be
  \begin{equation}
    \label{eq:10}
    \Sigma = \frac{1}{c\sinh c}
    \begin{pmatrix} f(s) & 0 \\ 0 & f(t) \end{pmatrix} \cdot
    \begin{pmatrix}
      \sinh cs \sinh c(1-s) &\sinh cs \sinh c(1-t)\\
      \sinh cs \sinh c(1-t) & \sinh ct \sinh c(1-t)
    \end{pmatrix}
    \cdot
    \begin{pmatrix} f(s) & 0 \\ 0 & f(t) \end{pmatrix}
  \end{equation}
  The matrix in the middle together with the factor $\frac{1}{c\sinh c}$ is
  the covariance matrix of the Ornstein-Uhlenbeck bridge on the interval
  $[0,1]$ with parameter $c$ (see Appendix).  Further computations
  reveal that
  \begin{equation*}
    \frac{1}{\sqrt{\det(\Sigma)}} = H^{(2)}_{\rho,c}(s,t,L(s),L(t))
  \end{equation*}
  since both sides equal
  \begin{equation*}
    \frac{c}{f(s)f(t)} \sqrt{\frac{\sinh{c}} {\sinh cs\sinh c(t-s)\sinh c(1-t)}}.
  \end{equation*}
  This completes the proof.
\end{proof}

We explain now that for $m$ points $0<t_1<t_2<\cdots<t_m<1$, the limit of the
corresponding $m$-dimensional marginal
$\paren[\big]{\tilde{X}^{(n)}_{t_1},\tilde{X}^{(n)}_{t_2},\ldots,\tilde{X}^{(n)}_{t_m}}$
is Gaussian, with covariance matrix $\Sigma=(\sigma_{i,j})$ defined by
\begin{equation*}
  \Sigma_{i,j} = \frac{1}{c\sinh{c}} f(t_i)f(t_j) \sinh ct_i 
  \sinh c(1-t_j) 
\end{equation*}
for all $i\leq j$.

A first approach would be to repeat similar computations as for the $2$-dimensional case, and observe that the matrix $\sigma$ is the inverse of the $m$-dimensional Hessian of $F^{(m)}_{\rho,c}(t_1,\dots,t_m)$.

Another way takes advantage of the Markov property satisfied by the process $(\tilde{X}_n)$. Denote by $P^{(n)}_{s,t}(x,y)$ the transition kernel of this Markov chain. 
The law of the $m$-tuple $\paren[\big]{\tilde{X}^{(n)}_{t_1},\tilde{X}^{(n)}_{t_2},\ldots,\tilde{X}^{(n)}_{t_m}}$ is given by

\begin{equation*}
P^{(n)}_{0,t_1}(0,x_1)P^{(n)}_{t_1,t_2}(x_1,x_2)\dots P^{(n)}_{t_{m-1},t_m}(x_{m-1},x_m).
\end{equation*}

Proposition \ref{prop:fluct2dim} implies the convergence of each $\sqrt{n} P^{(n)}_{t_i,t_{i+1}}(x_n,y_n)$ to the density of the kernel of the Ornstein-Uhlenbeck bridge as $n\rightarrow \infty$, and $(x_n)$, $(y_n)$ converge.
The convergence in law follows from Lebesgue dominated convergence theorem. Domination is ensured by Lemma \ref{lem:luf} below and unimodality stated in Lemma \ref{lem:unimodal}.  

\subsection{Tightness and proof of Theorem \ref{thm:2}}

We already proved the convergence of the finite-dimensional distributions. We need to show
that the sequence of the distribution is tight. For this we use Theorem 13.5 p.142 of
\cite{Billingsley}. The criterion is checked in Lemma \ref{lem:tightness}.

Before entering the proof, we need a geometric definition.

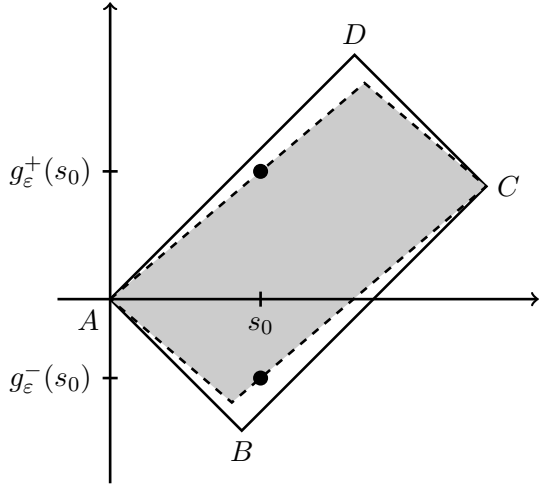
\begin{figure}[htb]
  \label{fig:parallelogram}
 \centering
  \begin{tikzpicture}[line width=1pt,scale=10]
    \renewcommand{\epsvalue}{.15} \renewcommand{\rhovalue}{.35}
   \path (0,0) coordinate (A) node[below left] {$A$};
    \path (.5*\rhovalue,-.5*\rhovalue) coordinate (B) node[below] {$B$};
    \path (.5,.5-\rhovalue) coordinate (C) node[right] {$C$};
    \path (.5-.5*\rhovalue,.5-.5*\rhovalue) coordinate (D) node[above] {$D$};
    \draw (A) -- (B) -- (C) -- (D) -- cycle;
    \path (intersection cs: first line={(0,0)--(1,1-\epsvalue)}, second
    line={(C)--(-.5,.5-\rhovalue-\epsvalue+1)}) coordinate
    (uppercorner); 
    \path (intersection cs: first line={(A)--(1,\epsvalue-1)},
    second line={(C)--(-.5,.5-\rhovalue+\epsvalue-1)}) coordinate
    (lowercorner); 
    \draw[dashed,fill=gray!40] (A) -- (uppercorner) -- (C) --
    (lowercorner) -- cycle;

    \newcommand{\snul}{.2}
    \draw (\snul,-.01) node[below] {$s_0$} -- (\snul,.01);
    \path (intersection cs: first line={(\snul,0)--(\snul,1)}, second
    line={(C)--(lowercorner)}) coordinate (gminus);
    \path (intersection cs: first line={(\snul,0)--(\snul,1)}, second
    line={(0,0)--(uppercorner)}) coordinate (gplus);
    \draw (gminus -| 0,0) ++(-.01,0) node[left] {$g_\epsilon^-(s_0)$} -- ++(.02,0);
    \draw (gplus -| 0,0) ++(-.01,0) node[left] {$g_\epsilon^+(s_0)$} -- ++(.02,0);
    \fill (gminus) circle (.01);
    \fill (gplus) circle (.01);

    \renewcommand{\axisplus}{0.07} 
    \draw[->] (0,-0.5*\rhovalue-\axisplus) -- (0,.5-.5*\rhovalue+\axisplus); 
    \draw[->] (-\axisplus,0) -- (.5+\axisplus,0);
  \end{tikzpicture}
  \caption{The $\epsilon$-parallelogram of a rectangle.}
\end{figure}

\begin{defi}
  Let $ABCD$ be a rectangle with sides having slopes $\pm1$ in the $(s,x)$
  coordinates.  For $\varepsilon\in(0,1)$, the $\varepsilon$-parallelogram of
  the rectangle $ABCD$ is the unique parallelogram with diagonal $AC$ and sides
  with slopes $1-\varepsilon$ and $-1+\varepsilon$ in the $(s,x)$
  coordinates. The $\varepsilon$-interior of the rectangle is the interior of
  the $\varepsilon$-parallelogram. The complement of the $\varepsilon$-interior
  in the rectangle is called the $\varepsilon$-boundary of the rectangle. See
  Figure \ref{fig:parallelogram}. For $s_0\in(0,1)$, the sides of the
  parallelogram intersect the straight line $s=s_0$ in two points. We denote the
  ordinates of these two points by $g_\varepsilon^+(s_0)$ and
  $g_\varepsilon^-(s_0)$.
\end{defi}

We first state the following useful fact.

\begin{fact}
\label{fact}
For all $\delta\in(0,1)$, for all $A>0$, there exists $\varepsilon$, such that for
  all rectangular boxes of side lengths $a$, $b$ satisfying $\rho=\frac{a}{a+b}
  \in(\delta,1-\delta)$, and all $c\in[-A,A]$, the limit shape $L_{\rho,c}$ is
  entirely included in the $\varepsilon$-interior of the box.
\end{fact}

The following lemma controls the function $F^{(1)}_{\rho,c}$ in an
$\varepsilon$-boundary of a rectangular box with side lengths $\rho$ and
$1-\rho$.

\begin{lem}\label{lem:F_boundary}
Let $\delta\in(0,1)$ and $A>0$. Take $\varepsilon$ as in Fact \ref{fact}.
  Then there exists a
  positive constant $C$ such that for all $\rho\in [\delta,1-\delta]$, and all $c\in[-A,A]$,  
  \begin{equation*}
    \forall (s,x)\in B_\rho^\varepsilon, \quad F^{(1)}_{\rho,c}(s,x) \leq -C s(1-s),
  \end{equation*}
  where $B_\rho$ is the (macroscopic) box with perimeter 1 and aspect ratio $\rho$.
\end{lem}
  
\begin{proof}
  For a fixed $s$, the function $x\mapsto F^{(1)}_{\rho,c}(s,x)$ is concave, and
  its maximum is reached at $x=L_{\rho,c}(s)$, and the point $(s,L_{\rho,c}(s))$
  is inside the $\varepsilon$-interior of the box. Therefore
  $F^{(1)}_{\rho,c}(s,x) \leq \max\{F^{(1)}_{\rho,c}(s,g_\varepsilon^+(s)),
  F^{(1)}_{\rho,c}(s,g_\varepsilon^-(s))\}$.
  

  In a neighborhood of $0$, $g_\varepsilon^+(s)=(1-\varepsilon)s$, and $s\mapsto
  F^{(1)}_{\rho,c}(s,(1-\varepsilon)s)\asymp s$ because it is differentiable and
  vanishes at $s=0$ with a non zero derivative.  The same argument applied to a neighborhood of $s=1$,
  and for $g_\varepsilon^-$ gives the result. The uniform bound on $C$ results from a compactness argument.
\end{proof}

\begin{lem}[L.U.F.]\label{lem:luf}
  Let $\delta\in(0,1)$ and $A>0$. Take $\varepsilon$ as in Fact \ref{fact}.
    For this $\varepsilon$, there exist two constants $M$ and $\kappa$ such that
  if $a+b=2n$ and $\frac{a}{a+b} \in(\delta,1-\delta)$ and for all $c\in[-A,A]$,
  then
  \begin{equation*}
    \forall\ 0\leq k \leq n,\quad \pr_{a,b}^{e^{-c/n}} (X_{2k} = 2\bigl\lfloor n L_{\rho,c}(\frac{k}{n})+y\sqrt{n}\bigr\rfloor) 
    \leq \frac{M}{\sqrt{n}}  \frac{\exp{\paren*{-\frac{y^2}{2\kappa s(1-s)}}}}{\sqrt{2\pi\kappa s(1-s)}}.
  \end{equation*}
  as soon as $(2k, \bigl\lfloor n
  L_{\rho,c}(\frac{k}{n})+y\sqrt{n}\bigr\rfloor)$ is in the
  $\varepsilon$-interior of the box.
\end{lem}

\begin{proof}
  Use Corollary~\ref{cor:1-marg}. Notice that $h_c(x,y) \asymp \sqrt \frac{x+y}{x y}$
  uniformly in $c$.
  Therefore
  \begin{equation*}
    H_{\rho,c}(s,x) \asymp \sqrt{\rho(1-\rho)}
    \sqrt{\frac{s}{(s-x)(s+x)}}\sqrt{\frac{1-s} 
      {(2\rho-s+x)(2-2\rho-s-x)}}.
  \end{equation*}
  The first term is bounded. The second term is of order $\frac{1}{\sqrt{s}}$ in a
  neighborhood of $(0,0)$ in the $\varepsilon$-interior of the box. A simple
  change of variable $(s,x)\mapsto(1-s,2\rho-1-x)$ exchanges the
  second and the third term, which shows that the third term is of order
  $\frac{1}{\sqrt{1-s}}$ in a neigborhood of $(1,1-\rho)$ in the
  $\varepsilon$-interior of the box. Therefore,
  \begin{equation*}
    H_{\rho,c}(s,x) \asymp \frac{1}{\sqrt{s(1-s)}},
  \end{equation*}
  as long as $(s,x)$ is in the $\varepsilon$-interior of the box.

  In order to bound the exponential term in Corollary~\ref{cor:1-marg}, we bound
  from below the absolute value of $\frac{\partial^2 F^{(1)}_{\rho,c}
    (s,x)}{\partial x ^2}$. Using that $S_c''(u)= \frac{c}{e^{c u}-1}$, we get
  \begin{equation*}
    \frac{\partial^2 F^{(1)}_{\rho,c} (s,x)}{\partial x ^2}= - 2c - \frac{c}{e^{c(s+x)}-1}
    -  \frac{c}{e^{c(s-x)}-1} -  \frac{c}{e^{c(2\rho-s+x)}-1} -
    \frac{c}{e^{c(2-2\rho-s-x)}-1} . 
  \end{equation*}
  For $s$ close to $0$, the main contribution comes from the second and the
  third term, that are both negative, and of order $s^{-1}$. For $s$ close to
  $1$, the main contribution comes from the fourth and the fifth term, that
  are both negative, and of order $(1-s)^{-1}$.
\end{proof}

\begin{lem} \label{lem:app_luf}
Let $\delta \in(0,1)$ and $A>0$, and $\varepsilon\in(0,1)$ given by Fact \ref{fact}. Then there exists a constant $C>0$ such that: for all $n\geq 1$, for all $c\in[-A,A]$, for all sequences of boxes $(B_n)$  with sides $a_n,b_n$ such that 
   $ a_n+b_n=2n, \frac{a_n}{2n}  \in(\delta,1-\delta),$
 \begin{equation*}
   \forall \lambda >0, \forall s\in[0,1],\quad \pr(\tilde{X}_s \not\in [-\lambda,
   \lambda] ; (X_{2ns},2ns) \in B_n^\varepsilon ) 
    \leq C\frac{(s(1-s))^2}{\lambda^4}
  \end{equation*}
\end{lem}

\begin{proof}
Lemma~\ref{lem:luf} 
gives the following bound: 
\begin{equation*}
\pr(\tilde{X}_s \not\in [-\lambda, \lambda] ; (X_{2ns},2ns) \in B_n^\varepsilon ) \leq \sum_{|y|\geq \lambda}  \frac{M}{\sqrt{n}}  \frac{\exp{(-\frac{y^2}{2\kappa s(1-s)}})}{\sqrt{2\pi\kappa s(1-s)}} ,
\end{equation*}
where the index of the sum $y$ runs on the set $\frac{1}{\sqrt{n}}\mathbb{Z}-\sqrt{n}L_{\rho,c}(s)$.

Comparing the sum on the right-hand side and the corresponding integral, we get
\begin{equation*}
\pr(\tilde{X}_s \not\in [-\lambda, \lambda] ; (X_{2ns},2ns) \in B_n^\varepsilon ) \leq M'  \int_{\mathbb{R}\setminus(-\lambda,\lambda)}\frac{\exp{(-\frac{y^2}{2\kappa s(1-s)}})}{\sqrt{2\pi\kappa s(1-s)}}\mathrm{d}s.
\end{equation*}
This last integral is equal to $\pr(|\mathcal{N} | \geq \frac{\lambda}{\sqrt{\kappa s
    (1-s)}})$, where $\mathcal{N}$ is a standard Gaussian variable. Conclude by using
Markov inequality for the fourth moment.
\end{proof}

We can now verify Billingsley's condition for tightness \cite{Billingsley}.

\begin{lem}\label{lem:tightness}
  Let $((a_n,b_n))_n$ be a sequence such that $a_n+b_n=2n$, and
  $\frac{a_n}{2n}=\rho +o(\frac{1}{\sqrt n})$. Then there exists a constant $C$
  such that for all $n>0$, for all $0\leq r \leq s \leq t\leq 1$, for
  all $\lambda>0$,
  \begin{equation}
    \pr_{a_n,b_n}^{e^{-c/n}} (|\tilde{X}_s -\tilde{X}_r|\geq \lambda ; |\tilde{X}_t -\tilde{X}_s|\geq \lambda) \leq \frac{C (t-r)^2}{\lambda^4}.
    \label{eq:tightness_bill}
  \end{equation}
\end{lem}

\begin{proof}
  Inequality \eqref{eq:tightness_bill} is automatically satisfied,  as soon as $|r-s|$ or
  $|t-s|$ is less than $1/n$, or $\lambda$ is greater than $\sqrt{n} (|r-s|\wedge |t-s|)
  $. We suppose now that none of these conditions are satisfied.

Introduce now the three following $\varepsilon$-interiors:
$B^\varepsilon$, the $\varepsilon$-interior of the box $B:=[0,a_n]\times [0,b_n]$, $B^\varepsilon_L$ the $\varepsilon$-interior of the box $B_L:=[0, ns-\frac{1}{2}X_{2sn}]\times [0,sn+\frac{1}{2}X_{2sn}]$ and $B^\varepsilon_R$ the $\varepsilon$-interior of the box $B_R:=[ ns-\frac{1}{2}X_{2sn},a_n]\times [sn+\frac{1}{2}X_{2sn},b_n]$. See Fig. \ref{fig:LRinteriors}.

\begin{figure}[htb]
  \centering
    \begin{tikzpicture}[line width=1pt,scale=13]
    \renewcommand{\epsvalue}{.3} \renewcommand{\rhovalue}{.4}
    \renewcommand{\axisplus}{0.03} 
    \coordinate (origin) at (0,0);
    \coordinate  (rightcorner) at (.5,.5-\rhovalue);
    \path (rightcorner) ++(-.1,.1*\epsvalue-.1) coordinate (epsrightlineup);
    \path (rightcorner) ++(-.1,-.1*\epsvalue+.1) coordinate (epsrightlinedown);
    \path (intersection cs: first line={(origin)--(1,1-\epsvalue)}, second
    line={(rightcorner)--(epsrightlinedown)}) coordinate (uppercorner); 
    \path (intersection cs: first line={(origin)--(1,\epsvalue-1)},
    second line={(rightcorner)--(epsrightlineup)}) coordinate (lowercorner); 
    \draw[dashed,fill=gray!40] (origin) -- (uppercorner) -- (rightcorner) --
    (lowercorner) -- cycle;
    \newcommand{\ns}{.28}\newcommand{\Xns}{.1}
    \newcommand{\etavalue}{.15}
    \path (\ns,\Xns) coordinate (center);
    \draw (\ns,-.01) node[below] {$ns$} -- (\ns,.01);
    \draw (-.01,\Xns) node[left] {$X_{ns}$} -- (.01,\Xns);
    \path (center) ++(.1,.1-.1*\etavalue) coordinate (centerlineup);
    \path (center) ++(.1,.1*\etavalue-.1) coordinate (centerlinedown);
    \path (rightcorner) ++(-.1,.1*\etavalue-.1) coordinate (etarightlineup);
    \path (rightcorner) ++(-.1,-.1*\etavalue+.1) coordinate (etarightlinedown);
    \path (intersection cs: first line={(origin)--(1,1-\etavalue)}, second
    line={(center)--(centerlinedown)}) coordinate (leftuppercorner); 
    \path (intersection cs: first line={(origin)--(1,\etavalue-1)},
    second line={(center)--(centerlineup)}) coordinate (leftlowercorner); 
    \draw[dashed,fill=blue!40] (origin) -- (leftuppercorner) -- (center) --
    (leftlowercorner) -- cycle;

    \path (intersection cs: first line={(center)--(centerlineup)}, second
    line={(rightcorner)--(etarightlinedown)}) coordinate (rightuppercorner);  
    \path (intersection cs: first line={(center)--(centerlinedown)}, second
    line={(rightcorner)--(etarightlineup)}) coordinate (rightlowercorner); 
    \draw[dashed,fill=red!40] (center) -- (rightuppercorner) -- (rightcorner) --
    (rightlowercorner) -- cycle;
    
    \draw[->] (0,-0.5*\rhovalue-\axisplus) -- (0,.5-.5*\rhovalue+\axisplus); 
    \draw[->] (-\axisplus,0) -- (.5+\axisplus,0);
    \draw[rotate=45] (origin) rectangle (rightcorner);
    \draw[rotate=45,blue] (origin) rectangle (center);
    \draw[rotate=45,red] (center) rectangle (rightcorner);
    \node at ($(origin)!0.5!(center)$) {$B_L^{\eta(\epsilon)}$};
    \node at ($(center)!0.5!(rightcorner)$) {$B_R^{\eta(\epsilon)}$};
  \end{tikzpicture}
  \caption{The $\epsilon$-interiors $B^{\varepsilon}$ (grey), $B_L^{\eta(\varepsilon)}$ (blue) and $B_R^{\eta(\varepsilon)}$ (red).}
  \label{fig:LRinteriors}
\end{figure}
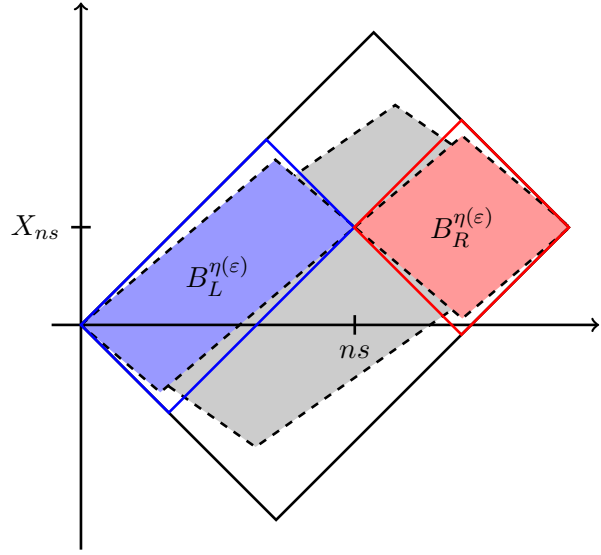

 
Let us first control the probability that the random interface exits the
$\varepsilon$-interior of the box. 
We use unimodality of the distribution of $X_{2ns}$, Corollary \ref{cor:1-marg}, and Lemma \ref{lem:F_boundary} to get:
 \begin{multline*}
 \pr((2ns,X_{2ns})\not\in B^\varepsilon) \leq 4\varepsilon n s(1-s) \pr( X_{2ns} = n g^\varepsilon_\pm(s)) \\
  \leq
 M \frac{ n s(1-s)}{\sqrt{n}} H_{\rho,c}(s,g^\varepsilon_{\pm}(s)) \exp (- n C s(1-s)). 
 \end{multline*}
 Using the uniform bound of $H_{\rho,c}\asymp \frac{1}{\sqrt{s(1-s)}}$ obtained in the beginning of the proof of Lemma \ref{lem:luf}, 
  we get
 \begin{equation*}
   \pr((2ns,X_{2ns})\not\in B^\varepsilon) \leq M \sqrt{n s (1-s)} \exp(-n C(s(1-s))) 
   \leq \frac{M}{(nC s(1-s))^2}.  
 \end{equation*}

 Using our assumption that $n \geq \frac{\lambda^2}{|s-r|^2\wedge|t-s|^2}$, we obtain:
 \begin{equation}
\pr((2ns,X_{2ns})\not\in B^\varepsilon) \leq  M    \frac{|r-s|^4\wedge|t-s|^4}{\lambda^4 s^2(1-s)^2}
   \leq M\frac{|r-s|^2\wedge|t-s|^2}{\lambda^4} \leq M\frac{|t-r|^2}{\lambda^4}.
   \label{eq:badbox}
   \end{equation}
   
   \begin{rem}
   The weaker condition $n\geq \frac{\lambda^2}{s^2 \wedge(1-s)^2}$ yields
   \begin{equation}
   \pr((2ns,X_{2ns})\not\in B^\varepsilon) \leq  M\frac{s^2\wedge(1-s)^2}{\lambda^4}
   \label{eq:badbox_gen}.
   \end{equation}
   \end{rem}
   
   Now we consider three cases depending on the values of $r$ and $t$.
   
   \textbf{Case 1:} $r < \frac{\lambda}{4\sqrt{n}}$.
   Then $|\tilde{X}_r| \leq 2r\sqrt{n} < \frac{\lambda}{2}$. Thus,
   \begin{multline*}
   \pr
   (|\tilde{X}_s -\tilde{X}_r|\geq \lambda ; |\tilde{X}_t -\tilde{X}_s|\geq \lambda) \leq \pr_{a_n,b_n}^{e^{-c/n}} (|\tilde{X}_s -\tilde{X}_r|\geq \lambda) \leq \pr
   (|\tilde{X}_s| \geq \frac{\lambda}{2}) \\
    \leq  \pr
    (|\tilde{X}_s| \geq \frac{\lambda}{2}, (ns,X_{ns})\in B^\varepsilon ) + \pr
    ((ns,X_{ns})\not\in B^\varepsilon ).
   \end{multline*}
   
   But we just saw in \eqref{eq:badbox} that
    \begin{equation*}
   \pr
   ((2ns,X_{2ns})\not\in B^\varepsilon ) \leq M\frac{|t-r|^2}{\lambda^4}.
   \end{equation*}
   Moreover, according to Lemma \ref{lem:app_luf}, 
   \begin{equation*}
   \pr
   (|\tilde{X}_s| \geq \frac{\lambda}{2}, (2ns,X_{2ns})\in B^\varepsilon ) \leq C\frac{s^2(1-s)^2}{\lambda^4} \leq \frac{C s^2}{\lambda^4}.
   \end{equation*}
   Since $r < \frac{\lambda}{4\sqrt{n}} < \frac{s}{4}$, we get that $s < \frac{4}{3} (s-r)$, and thus
   \begin{equation*}
   \pr
   (|\tilde{X}_s| \geq \frac{\lambda}{2}, (2ns,X_{2ns})\in B^\varepsilon ) \leq \frac{C (s-r)^2}{\lambda^4}\leq  \frac{C (t-r)^2}{\lambda^4}.
   \end{equation*}
   And the tightness condition \eqref{eq:tightness_bill} is satisfied in that case.
   
 \textbf{Case 2:} $1-t < \frac{\lambda}{4\sqrt{n}}$. This case is treated similarly as the previous one, and \eqref{eq:tightness_bill} is again satisfied.
 
 \textbf{Case 3:} $r \geq \frac{\lambda}{4\sqrt{n}}$ and $1-t \geq \frac{\lambda}{4\sqrt{n}}$. This is the generic situation.

 Conditional on $\{(2ns,X_{2ns}) \in B^\varepsilon\}$, the boxes $B_L$ and $B_R$ have an
 aspect ratio in the interval $(\varepsilon,1-\varepsilon)$. Applying Lemma
 \ref{lem:app_luf} with $\delta=\varepsilon$ in the subboxes $B_L$ and $B_R$, we get the
 existence of an $\eta(\varepsilon)$ such that with sufficient probability $(2nr,X_{2nr})$
 (resp. $(2nt,X_{2nt})$) is in $B_L^{\eta(\varepsilon)}$
 (resp. $B_R^{\eta(\varepsilon)}$). More precisely, the assumptions on $r$ and $t$ allow to repeat the argument to derive \eqref{eq:badbox_gen} with the proper scaling in the subboxes $B_L$ and $B_R$ and get for some constant $M$:
 \begin{align}
   \pr
   ( (2nr,X_{2nr})\not\in B^{\eta(\varepsilon)}_L |(2ns,X_{2ns}) \in B^\varepsilon ) &\leq M \frac{
     (\frac{r}{s})^2 \wedge (\frac{s-r}{s})^2}{(\lambda/\sqrt{s})^4} \leq M
   \frac{(s-r)^2}{\lambda ^4} \label{eq:badboxL},\\ 
   \pr
   ( (2nt,X_{2nt})\not\in B^{\eta(\varepsilon)}_R |(2ns,X_{2ns}) \in B^\varepsilon ) &\leq
   M \frac{ (\frac{t-s}{1-s})^2 \wedge
     (\frac{1-t}{1-s})^2}{(\lambda/\sqrt{1-s})^4} \leq M
   \frac{(t-s)^2}{\lambda ^4}\label{eq:badboxR}.
 \end{align}

By the Markov property of $X$, the variables $X_r$ and $X_t$
are independent conditional on the value of $X_s$. Hence we can write
\begin{multline}
  \pr
  (|\tilde{X}_s -\tilde{X}_r|\geq \lambda ; |\tilde{X}_t -\tilde{X}_s|\geq \lambda)= \\
  = \sum_{\substack{j,k,l\\ |j-k|\geq \lambda\\ |l-k|\geq \lambda}}
  \pr
  (\tilde{X}_s = k)\pr
  (\tilde{X}_r=j | \tilde{X}_s = k) \pr
  (\tilde{X}_t=l |\tilde{X}_s = k),
  \label{eq:tight_split}
\end{multline}
where $j,k,l$ runs through all possible values for $\tilde{X}_r,\tilde{X}_s,
\tilde{X}_t$ (two successive values of $j,k,l$ differ by $1/\sqrt{n}$).
 
Let $E_{r,s,t}^\epsilon = \{(2ns,X_{2ns}) \in B^\varepsilon\} \cap
\{(2nr,X_{2nr})\in B^{\eta(\varepsilon)}_L \}\cap \{(2nt,X_{2nt})\in
B^{\eta(\varepsilon)}_R\}$. From Inequalities \eqref{eq:badbox}, \eqref{eq:badboxL},
\eqref{eq:badboxR}, there exists a constant $M$ such that:
\begin{equation*}
  1-\pr (E_{r,s,t}^\epsilon) \leq  M \frac{(t-r)^2}{\lambda^4}.
\end{equation*}
As a consequence, it suffices to bound the sum in \eqref{eq:tight_split} for
values of $(j,k,l)$ such that
\begin{equation*}
\{\tilde{X}_r =j,\tilde{X}_s=k, \tilde{X}_t=l\} \subset 
 \{|\tilde{X}_s -\tilde{X}_r|\geq \lambda\}\cap\{ |\tilde{X}_t -\tilde{X}_s|\geq
 \lambda\}\cap E_{r,s,t}^\varepsilon.
 \end{equation*}



 By a scaling argument, when looking at what happens on the left of $s$,
\begin{equation*}
  \pr_{a_n,b_n}^{e^{-c/n}} (\tilde{X}_r = j | X_{ns}) = 
  \pr_{\frac{2ns-X_{2ns}}{2},\frac{2ns+X_{2ns}}{2}}^{e^{-c/n}}\left(\tilde{X}_{\frac{r}{s}} 
= \frac{j}{\sqrt{s}}+\sqrt{\frac{n}{s}}L_{\rho,c}(r)-\sqrt{ns}L_{\frac{2ns-X_{2ns}}{4ns},sc}(\frac{r}{s})
\right).
\end{equation*}
Note that on the right hand side, $\tilde{X}_{r/s}$ is defined with respect to
the limit shape inside the box $[0,\frac{2ns-X_{2ns}}{4ns}]\times[0,\frac{2ns+X_{2ns}}{4ns}]$. See
Fig. \ref{fig:sticking}.

\begin{figure}[htb]
\centering
\begin{tikzpicture}[line width=1pt,scale=12]
  \renewcommand{\cvalue}{2}  \renewcommand{\rhovalue}{.6}
  \renewcommand{\axisplus}{0.03} 
  \newcommand{\svalue}{.28}\newcommand{\Xsvalue}{-0.04}
  \newcommand{\rvalue}{.15}\newcommand{\tvalue}{.42}
  \path (\svalue,\Xsvalue) coordinate (center);
  \path (.5,.5-\rhovalue) coordinate (rightcorner);
  \draw (\svalue,.01) node[above] {$s$} -- (\svalue,-.01);
  \draw (\rvalue,.01) node[above] {$r$} -- (\rvalue,-.01);
  \draw (\tvalue,.01) node[above] {$t$} -- (\tvalue,-.01);
  \draw (.5,.01) node[above] {\small $1$} -- (.5,-.01);
  \draw[rotate=45,dashed] (origin) rectangle (center);
  \draw[rotate=45,dashed] (center) rectangle (rightcorner);
  \draw[->] (0,-0.5*\rhovalue-\axisplus) -- (0,.5-.5*\rhovalue+\axisplus); 
  \draw[->] (-\axisplus,0) -- (.5+\axisplus,0);
  \path (0,0) coordinate (origin);
  \draw[rotate=45] (origin) rectangle (rightcorner);

  \path[rotate=-45] (center |- 0,0) node[above=5] {$B_L$};
  \path[rotate=-45] (center -| rightcorner) node[above=5] {$B_R$};

  \begin{scope}[smooth,color=purple,domain=0:0.5,samples=50]
    \draw plot (\x,{\limitshape{\cvalue}{\rhovalue}{\x}});
    \draw[scale={2*\svalue}] plot
    (\x,{\limitshape{2*\svalue*\cvalue}{0.5*(1-(\Xsvalue/\svalue))}{\x}});
    \draw[shift={(\svalue,\Xsvalue)},scale=2*(0.5-\svalue)] plot
    (\x,{\limitshape{2*\cvalue*(0.5-\svalue)}{0.5*(1-((0.5-\rhovalue-\Xsvalue)/
        (0.5-\svalue)))}{\x}}); 
  \end{scope}

  \pgfmathsetmacro{\upperend}
  {2*\svalue*\limitshape{2*\svalue*\cvalue}{.5*(1-(\Xsvalue/\svalue))}{\rvalue/2/\svalue}}
  \pgfmathsetmacro{\lowerend}
  {\limitshape{\cvalue}{\rhovalue}{\rvalue}};
  \draw[<->] (\rvalue,\upperend) -- (\rvalue,\lowerend);
\end{tikzpicture}
\caption{    Limit shapes in the box $B$ and the subboxes $B_L$ and
  $B_R$. The double arrow represents the difference of the limit shapes in $B$ and $B_L$
  at time $r$, whose expression, $2 n L_{\rho,c} (r) - 2 n s
  L_{\frac{2ns-X_{2ns}}{4ns},sc}\paren*{\frac{r}{s}} $, appears in the sticking condition
  \eqref{eq:stickingL}.}
\label{fig:sticking} 
\end{figure}
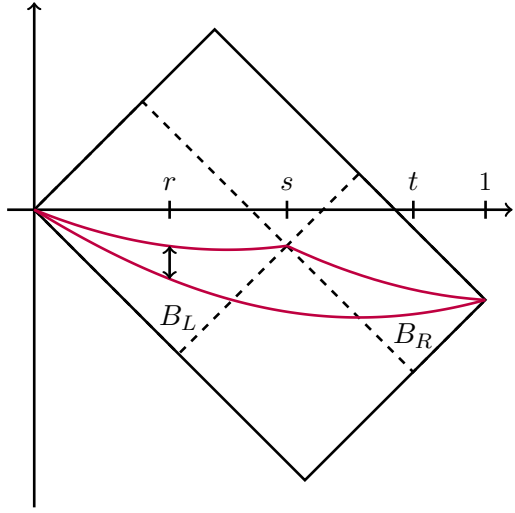

When summing over the values of $j$ such that $|j-\tilde{X}_s|\geq \lambda$, one can get an upper bound for $\pr_{a_n,b_n} (|\tilde{X}_r - \tilde{X}_s|\geq \lambda | \tilde{X}_s)$ using Lemma  \ref{lem:app_luf} as long as the limit shape in the original box $B$ and the one in $B_L$ are sufficiently close to each other. We express this proximity by the following sticking condition for $B_L$:
\begin{equation}
\frac{1}{ \sqrt{s}}[\tilde{X}_s-\lambda,\tilde{X}_s+\lambda] 
+\sqrt{\frac{n}{s}} L_{\rho,c} (r) - \sqrt{n s} L_{\frac{2ns-X_{2ns}}{2ns},sc}(\frac{r}{s}) 
\supset \left[-\frac{\lambda}{\sqrt{2s}},\frac{\lambda}{\sqrt{2s}}\right].
\label{eq:stickingL}
\end{equation}

Indeed, when the sticking condition \eqref{eq:stickingL} is satisfied, $\pr_{a_n,b_n}^{e^{-c/n}} (|\tilde{X}_r - \tilde{X}_s|\geq \lambda | \tilde{X}_s )  $ takes the form:
\begin{multline*}
  \pr_{\frac{2ns-X_{2ns}}{2},\frac{2ns+X_{2ns}}{2}}^{e^{-c/n}}\left(\tilde{X}_{\frac{r}{s}} \not\in \frac{1}{ \sqrt{2s}}[\tilde{X}_s-\lambda,\tilde{X}_s+\lambda] 
+\sqrt{\frac{n}{s}} L_{\rho,c} (r) - \sqrt{n s} L_{\frac{2ns-X_{2ns}}{4ns},sc}(\frac{r}{s}) \right)\\
\leq \pr_{\frac{2ns-X_{2ns}}{2},\frac{2ns+X_{2ns}}{2}}^{e^{-c/n}}\left(\tilde{X}_{\frac{r}{s}} \not\in \left[-\frac{\lambda}{\sqrt{2s}},\frac{\lambda}{\sqrt{2s}}\right] \right)
\leq M \frac{(\frac{r}{s} \frac{s-r}{s})^2}{(\lambda/\sqrt{2s})^4} \leq  4 M \frac{(t-r)^2}{\lambda^4}  
\end{multline*}
by Lemma \ref{lem:app_luf} in the small box.

We follow the same idea for the right box $B_R$. 
Whenever the sticking condition for $B_R$ is satisfied, 
\begin{equation*}
\pr_{a_n,b_n}^{e^{-c/n}} (|\tilde{X}_t - \tilde{X}_s|\geq \lambda | \tilde{X}_s )  \leq 4 M \frac{(t-r)^2}{\lambda^4}.
\end{equation*}

Therefore, conditional on one of the sticking conditions to be satisfied, the probability of $\{|\tilde{X}_s -\tilde{X}_r|\geq \lambda\}\cap\{ |\tilde{X}_t -\tilde{X}_s|\geq \lambda\}$ is bounded by $4 M\frac{(t-r)^2}{\lambda^4}$.

Let us now investigate the probability that both of the sticking conditions fail. It is bounded by the probability that the one for $B_L$ fails.

Using Equation \eqref{eq:27} of Lemma \ref{lem:1}, we replace $\frac{1}{s} L_{\rho,c}(r)$ in \eqref{eq:stickingL} by
 $ L_{\frac{s-L_{\rho,c}(s)}{2s},sc}\left(\frac{r}{s}\right)$.

Condition \eqref{eq:stickingL} is equivalent to
\begin{equation*}
\left| \tilde{X}_s + s\sqrt{n} \left( L_{\frac{s-L_{\rho,c}(s)}{2s},sc}\left(\frac{r}{s}\right)- L_{\frac{s-L_{\rho,c}(s)-\tilde{X}_s n^{-1/2}}{2s},sc}\left(\frac{r}{s}\right) \right) \right| \leq \frac{\lambda}{2}.
\end{equation*}

By the mean value theorem applied to $\rho\mapsto L_{\rho,sc}(r/s)$, there exists a $\bar\rho $ in the interval $(\frac{s-L_{\rho,c}(s)}{2s},\frac{s-L_{\rho,c}(s)-\tilde{X}_s n^{-1/2}}{2s})$ such that
\begin{equation*}
s\sqrt{n} \left( L_{\frac{s-L_{\rho,c}(s)}{2s},sc}\left(\frac{r}{s}\right)- L_{\frac{s-L_{\rho,c}(s)-\tilde{X}_s n^{-1/2}}{2s},sc}\left(\frac{r}{s}\right)\right) = 
\frac{1}{2}\tilde{X}_s \left.\frac{\partial L_{\rho,cs}(r/s)}{\partial \rho}\right|_{\rho=\bar\rho}.
\end{equation*}

But differentiating \eqref{eq:7} with respect to $\rho$, for generic values of $\rho,c,t$ yields
\begin{multline*}
  0\leq 1+\frac{1}{2}\frac{\partial L_{\rho,c}(t)}{\partial \rho} = \frac{e^{c(2\rho-1)} \sinh(c(1-t))}{\sinh(ct) +e^{c(2\rho-1)} \sinh(c(1-t))}\\
\leq 1 \wedge e^{c(2\rho-1)} \frac{\sinh(c(1-t))}{\sinh(ct)}
\leq 1 \wedge K \frac{1-t}{t}
\end{multline*}
with a constant $K$ which works for all $c\in[-A,A]$, all $\rho\in(0,1)$, and all $t\in[0,1]$. 

As a consequence, the sticking condition \eqref{eq:stickingL} is satisfied as soon as

\begin{equation*}
|\tilde{X}_s| (1\wedge K \frac{s-r}{r}) \leq \frac{\lambda}{2}.
\end{equation*}

Therefore the probability that \eqref{eq:stickingL} is not verified is less than
\begin{equation*}
P(|\tilde{X}_s| (1\wedge K \frac{s-r}{r}) > \frac{\lambda}{2}, (2sn,X_{2sn})\in B^\varepsilon),
\end{equation*}
which, by Lemma \ref{lem:app_luf} is bounded by $\lambda^{-4} s^2(1-s)^2 (1\wedge K^4(\frac{s-r}{r})^4)$. This bound is less than some constant times $\lambda^{-4}(s-r)^2 \leq \lambda^{-4}(t-r)^2$, as one can check in both regimes $\bigl(\frac{s-r}{r}\bigr) \lessgtr K$.
This finishes the proof of Lemma~\ref{lem:tightness}.
\end{proof}

\section{The unbounded case}

In this section, we relax the constraint of remaining in a box. For any $q\in(0,1)$, we define a probability measure $\pr^q$ on \emph{all} partitions, by
\begin{equation*}
  \pr^q(\lambda)= \frac{1}{Z(q)} q^{|\lambda|}, \text{for all partitions $\lambda$,}
\end{equation*}
where $Z(q) =\prod_{i=1}^{\infty} (1-q^i)^{-1}$ is the generating function of all partitions.

When $q$ goes to 1, a law of large number for the shape of the random partition occurs. Namely, when rescaled by $(1-q)$, the boundary of $\lambda$ converges to a deterministic curve.
\begin{thm}[Vershik]\label{thm:vershik_shape}
  Let $L_\infty(s)= \log(2\cosh s)$,  and $(X_s)_{s\in\mathbb{R}}$ be the piecewise linear function describing the boundary of $\lambda$ in the $(s,x)$-coordinates.
  Then, for all $\varepsilon>0$,
  \begin{equation*}
    \lim_{q\to 1} P_q(\sup_{s\in\mathbb{R}} |(1-q)X_{s(1-q)^{-1}} - L_\infty(s)| > \varepsilon) = 0.
  \end{equation*}
\end{thm}

\subsection{Link between $L_\infty$ and $L_{\rho,c}$}

Petrov recently discussed in \cite{Petrov} the link between the limit shapes  obtained above and the infinite shape of the unbounded problem. We explain it in this short section. 

The family of limit shapes obtained above has the following property: if we fix
$s_0\in (0,1)$ and take the point $(s_0,L_{\rho,c}(s_0))$ as the
right corner of a new family of bounding boxes, then the limit shape of this new
problem is simply the restriction of $L_{\rho,c}$ to the interval $[0,s_0]$,
rescaled by $s_0$ (See Lemma \ref{lem:1}). 

\begin{figure}[htb]
  \centering
  \begin{tikzpicture}[line width=1pt,scale=10]
    \renewcommand{\cvalue}{3} \renewcommand{\rhovalue}{.45}
    \draw[smooth,color=purple] plot[domain=0:0.5,samples=50] 
    (\x,{\limitshape{\cvalue}{\rhovalue}{\x}});
    \draw (0,0) -- (\rhovalue/2,-\rhovalue/2) -- (0.5,0.5-\rhovalue) --
    (1/2-\rhovalue/2,1/2-\rhovalue/2) -- cycle; 
    \pgfmathsetmacro{\ycoord}{\limitshape{\cvalue}{\rhovalue}{0.35}}
    \path (0.35,\ycoord) coordinate (P);
    \draw[rotate=45,dashed] (0,0) rectangle (P);
  \end{tikzpicture}
  \caption{The limit shape in a subbox.}
  \label{fig:subbox}
\end{figure}
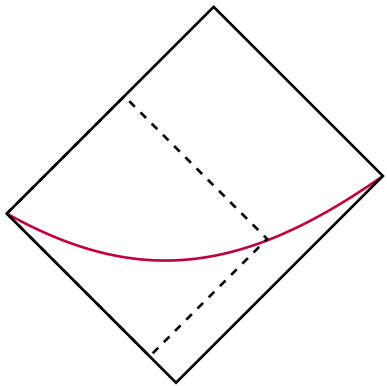

We are therefore led to the idea that there exists some
`inverse limit' of these
curves, from which they can all be recovered by restriction. This universal curve is  precisely $L_\infty$ from Theorem \ref{thm:vershik_shape}.
Indeed, for every $\rho\in(0,1)$ and $c\in\mathbb{R}_+^*$, one can find two real numbers $s_0 < s_1$ such that the restriction of $L_\infty$ to $[s_0,s_1]$ is (up to an affine transformation) $L_{\rho,c}$ (see Figure \ref{fig:restriction}):
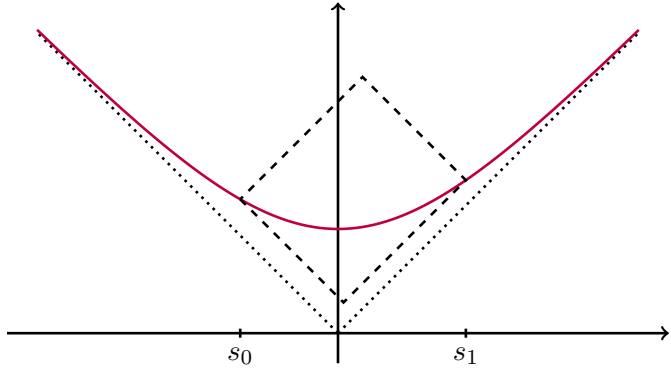
\begin{figure}[h]
  \centering
  \begin{tikzpicture}[line width=1pt,scale=2]
    \renewcommand{\axisplus}{0.2}
    \draw[->] (0,-\axisplus) -- (0,2+\axisplus); \draw[->] (-2-\axisplus,0) --
    (2+\axisplus,0); 
    \draw[dotted] (0,0) -- (-2,2); \draw[dotted] (0,0) -- (2,2); 
    \renewcommand{\cvalue}{1} 
    \newcommand{\inflimshape}[2]{(1/(#1)*ln(exp(#1*#2)+exp(-#1*#2)))}
    \draw[smooth,color=purple] plot[domain=-2:2,samples=50]
    (\x,{\inflimshape{\cvalue}{\x}});
    \newcommand{\leftx}{-0.65} \newcommand{\rightx}{.85}
    \pgfmathsetmacro{\lefty}{\inflimshape{\cvalue}{\leftx}}
    \pgfmathsetmacro{\righty}{\inflimshape{\cvalue}{\rightx}}
    \path (\leftx,\lefty) coordinate (leftcorner); \path
    (\rightx,\righty) coordinate (rightcorner);
    \begin{scope}[rotate=45]
      \draw[dashed] (leftcorner) rectangle (rightcorner);
    \end{scope}
    \draw (\leftx,-.03) -- (\leftx,.03); \path (\leftx,-.15) node {$s_0$}; \draw
    (\rightx,-.03) -- (\rightx,.03); \path (\rightx,-.15) node {$s_1$};
  \end{tikzpicture}
  \caption{Limit shapes are restrictions of the universal curve.}
  \label{fig:restriction}
\end{figure}
\begin{equation}
\forall t\in[0,1], \quad L_{\rho,c}(t) = \frac{L_{\infty}(s_0+t(s_1-s_0))-L_{\infty}(s_0)}{s_1-s_0}.
\end{equation}

For given values of $s_0 <s_1$, the parameter $\rho$  is given by
\begin{equation*}
  \rho=\frac{1}{2}-\frac{L_\infty(s_1)-L_\infty(s_0)}{2(s_1-s_0)},
\end{equation*}
and $c$ is found for example by comparing the slopes of the limit shape in the corner and identifying $\left.\frac{\partial L_\infty}{\partial s}\right|_{s=s_0}(s_1-s_0)$ and $\left.\frac{\partial L_{\rho,c}}{\partial s} \right|_{s=0} = \frac{e^{c(2\rho-1)}-\cosh c}{\sinh c}$, which is an increasing function of $c$.

The limit case $c\to 0$ correspond to $s_1\to s_0$. For negative values of $c$ one can take advantage of the symmetry of the model.

\subsection{Fluctuations around $L_\infty$}

As in the boxed case, we can compute explicitely the two point distribution. The analogue of Proposition \ref{prop:2marg} in the unbounded case is
\begin{prop}
  \label{prop:2marg_unbounded}
  The 2-dimensional marginal of $(X_{2k})_{ k\in\mathbb{Z}}$ is given by
  \begin{multline*}
    \pr^q(X_{2k}=2i,X_{2l}=2j)=\\
    \frac{Z_{k+i,\infty}(q) Z_{j+l-i-k,l+i-k-j}(q)Z_{j-l,\infty}(q)}{Z(q)}
    q^{(i-k)(k+i)+(l+j-k-i)(j-l)}.
  \end{multline*}
  where $Z_{a,\infty}(q)=\lim_{b\to \infty} Z_{a,b}(q) = \prod_{j=1}^a (1-q^j)^{-1}$.

\end{prop}

The asymptotic analysis of this formula as $q\to 1$ relies on q-Stirling's Formula (Corollary \ref{prop:qstirling}) and goes as in Section \ref{subsec:fluc2pt}.

One can then deduce the following theorem:

\begin{thm}\label{thm:fluct_infty}
  As $q$ goes to $1^-$, the random function
  \begin{equation*}
    s\mapsto \frac{\sqrt{2}\cosh s}{\sqrt{1-q}} \left\{(1-q)X_{s(1-q)^{-1}}-L_{\infty}(s)\right\}
  \end{equation*} converges weakly in $D$ to the
  two-sided stationary Ornstein-Uhlenbeck process $(Y_s)_{s\in \mathbb{R}}$, which is the Gaussian process on $\mathbb{R}$ 
  with covariance 
\begin{equation*}
  \forall(s,t)\in\mathbb{R}^2, \quad  \E\brac{Y_s Y_t} =  e^{-|t-s|}.
\end{equation*}

\end{thm}

When conditioning the process $Y_s / \sqrt{2}\cosh s$ to have zero integral over $\mathbb{R}$, one obtains a new centered Gaussian process whose covariance is equal to
\begin{equation*}
  \frac{e^{-|t-s|}}{2\cosh s \cosh t} - \frac{\pi^2}{6}h(s) h(t),
\end{equation*}
where $h(t)=\frac{6}{\pi^2}t \tanh t - \log(2\cosh t)$.

Expressed in the original coordinates (by rotating back the picture by 45 degrees), this covariance is the one obtained by Pittel \cite{Pittel} who deals with the microcanonical ensemble of partitions (with a fixed area $|\lambda|=n\to\infty$). The fluctuations in that case are thus given by this Gaussian, but non-Markov process.

%
%
%
%

\section{Appendix: The Ornstein-Uhlenbeck Bridge}
\label{sec:append-ornst-uhlenb}

In this short appendix, we give a description of the Ornstein-Uhlenbeck
bridge. Let $(B_t)$ be a Brownian motion on $[0,\infty)$. The Ornstein-Uhlenbeck
process $(Z_t)$ is the Gaussian random process defined by 
\begin{equation}
  \label{eq:55}
  Z_t = B_t - c \int_0^t Z_s \d{s}=e^{-ct}\int_0^t e^{cs}\d{B_s}.
\end{equation}
Using the right-hand side of \eqref{eq:55}, we can also represent $Z_t$ as: 
\begin{equation}
  \label{eq:60}
  Z_t = e^{-ct}B_{\int_0^t e^{cs}\d{s}} = e^{-ct}B_{\frac{e^{2ct}-1}{2c}}.
\end{equation}
From \eqref{eq:60} we derive the covariance. For $0\leq s\leq t$,
\begin{equation}
  \label{eq:61}
  \E\brac*{Z_s Z_t} = e^{-c(s+t)} \frac{e^{2cs}-1}{2c} =
  e^{-ct}\frac{\sinh(cs)}{c}. 
\end{equation}
`Tying down' this process at $t=1$ gives the Ornstein-Uhlenbeck $(Y_t)_{t\in[0,1]}$
bridge of length $1$. Set
\begin{equation}
  \label{eq:32}
  h(t) = \frac{\sinh(ct)}{\sinh{c}}
\end{equation}
and define
\begin{equation}
  \label{eq:31}
  Y_t = Z_t - h(t) Z_{1}.
\end{equation}
Observe that, for $0\leq t\leq 1$,
\begin{align*}
  \E\brac*{Y_t Z_{1}} &= \E\brac*{Z_t Z_{1}-h(t) Z_{1}^2} \\ &=
  e^{-c}\frac{\sinh(ct)}{c} - \frac{\sinh(ct)}{\sinh{c}}
  e^{-c}\frac{\sinh{c}}{c}\\ &=0,
\end{align*}
showing that $Y_t$ is independent of $Z_{1}$. Now let $f$ be any test function
and consider
\begin{align*}
  \label{eq:33}
  \E\brac*{f(X_t)\mid X_{1}=0} &= \E\brac*{f(Y_t+h(t)X_{1}) \mid
    X_{1}=0} \\ &= \E\brac*{f(Y_t)\mid X_{1}=0} \\ &= \E\brac{f(Y_t)},
\end{align*}
where the last equality is a consequence of the independence just proved. Thus,
the process $(Y_t)$ has the probability law of the Ornstein-Uhlenbeck process
conditional on $X_{1}=0$. Its covariance is
\begin{align*}
  \E\brac{Y_sY_t} &= \E\brac*{Y_s(Z_t-h(t) Z_{1})} = \E\brac{Y_s Z_t} \\ &=
  \E\brac{Z_s Z_t} -h(s)\E\brac{Z_{1} Z_t} =
  \frac{\sinh(cs)}{c\sinh{c})}
  \paren*{e^{-ct}\sinh{c}-e^{-c}\sinh(ct)} \\ 
  &= \frac{\sinh(cs)\sinh(c(1-t))}{c\sinh{c}},
\end{align*}
cf. the middle matrix in \eqref{eq:10}.

Notice that this covariance is also the Green function of a Brownian motion on $[0,1]$ killed
an exponential rate $c^2$ \cite{Chan_et_al}.
Let us mention also, that this process was used by C. Donati in her solution of Buffon-Synge problem concerning the typical distance between the extremities of a string of given length, thrown at random \cite{Donati-Martin}. 
\vspace{1cm}

\textbf{Acknowledgements. } The authors are pleased to thank Nikolai Reshetikhin for having proposed the question solved in this paper, Marc Yor for having recognized the covariance we obtained as that of the Ornstein-Uhlenbeck bridge, and Richard Kenyon for showing us the way to get the limit shapes we found from the infinite limit shape. The work of the third author was partially supported by the ANR  MEMEMO 
grant.


\bibliographystyle{plain}
\bibliography{partitions-rectangle}	

\end{document}